\documentclass[11pt]{article}

\begin{document}

%%-I include macros files here-%%

%%dung le's macros

\newcommand{\newc}{\newcommand}

%%index

\renewcommand{\theequation}{\thesection.\arabic{equation}}
\newc{\eqnoset}{\setcounter{equation}{0}}
\newcommand{\myref}[2]{#1~\ref{#2}}

\newcommand{\mref}[1]{(\ref{#1})}
\newcommand{\reflemm}[1]{Lemma~\ref{#1}}
\newcommand{\refrem}[1]{Remark~\ref{#1}}
\newcommand{\reftheo}[1]{Theorem~\ref{#1}}
\newcommand{\refdef}[1]{Definition~\ref{#1}}
\newcommand{\refcoro}[1]{Corollary~\ref{#1}}
\newcommand{\refprop}[1]{Proposition~\ref{#1}}
\newcommand{\refsec}[1]{Section~\ref{#1}}
\newcommand{\refchap}[1]{Chapter~\ref{#1}}

%%environments
\newcommand{\beq}{\begin{equation}}
\newcommand{\eeq}{\end{equation}}
\newcommand{\beqno}[1]{\begin{equation}\label{#1}}

\newcommand{\barr}{\begin{array}}
\newcommand{\earr}{\end{array}}

\newc{\bearr}{\begin{eqnarray*}}
\newc{\eearr}{\end{eqnarray*}}

\newc{\bearrno}[1]{\begin{eqnarray}\label{#1}}
\newc{\eearrno}{\end{eqnarray}}

\newc{\non}{\nonumber}
\newc{\nol}{\nonumber\nl}

\newcommand{\bdes}{\begin{description}}
\newcommand{\edes}{\end{description}}
\newc{\benu}{\begin{enumerate}}
\newc{\eenu}{\end{enumerate}}
\newc{\btab}{\begin{tabular}}
\newc{\etab}{\end{tabular}}

%%\newtheorem{theorem}{Theorem}
%%\newtheorem{defi}[theorem]{Definition}
%%\newtheorem{lemma}{Lemma}[section]
%%\newtheorem{rem}[lemma]{Remark}
%%\newtheorem{exam}[lemma]{Example}
%%\newtheorem{propo}[theorem]{Proposition}
%%\newtheorem{corol}[theorem]{Corollary}

%%Unmark these for uniform indexing

\newtheorem{theorem}{Theorem}[section]
\newtheorem{defi}[theorem]{Definition}
\newtheorem{lemma}[theorem]{Lemma}
\newtheorem{rem}[theorem]{Remark}
\newtheorem{exam}[theorem]{Example}
\newtheorem{propo}[theorem]{Proposition}
\newtheorem{corol}[theorem]{Corollary}

\renewcommand{\thelemma}{\thesection.\arabic{lemma}}

\newcommand{\btheo}[1]{\begin{theorem}\label{#1}}
\newc{\brem}[1]{\begin{rem}\label{#1}\em}
\newc{\bexam}[1]{\begin{exam}\label{#1}\em}
\newc{\bdefi}[1]{\begin{defi}\label{#1}}
\newcommand{\blemm}[1]{\begin{lemma}\label{#1}}
\newcommand{\bprop}[1]{\begin{propo}\label{#1}}
\newcommand{\bcoro}[1]{\begin{corol}\label{#1}}
\newcommand{\etheo}{\end{theorem}}
\newcommand{\elemm}{\end{lemma}}
\newcommand{\eprop}{\end{propo}}
\newcommand{\ecoro}{\end{corol}}
\newc{\erem}{\end{rem}}
\newc{\eexam}{\end{exam}}
\newc{\edefi}{\end{defi}}

\newc{\rmk}[1]{{\bf REMARK #1: }}
\newc{\DN}[1]{{\bf DEFINITION #1: }}

\newcommand{\bproof}{{\bf Proof:~~}}
\newc{\eproof}{{\vrule height8pt width5pt depth0pt}\vspace{3mm}}

%%symbols

\newcommand{\rarr}{\rightarrow}
\newcommand{\Rarr}{\Rightarrow}
\newcommand{\tru}{\backslash}
\newc{\bfrac}[2]{\dspl{\frac{#1}{#2}}}

%%space

\newc{\nl}{\vspace{2mm}\\}
\newc{\nid}{\noindent}

%%operations

\newcommand{\oneon}[1]{\frac{1}{#1}}
\newcommand{\dspl}{\displaystyle}
\newc{\grad}{\nabla}
\newc{\Div}{\mbox{div}}
\newc{\pdt}[1]{\dspl{\frac{\partial{#1}}{\partial t}}}
\newc{\pdn}[1]{\dspl{\frac{\partial{#1}}{\partial \nu}}}
\newc{\pdNi}[1]{\dspl{\frac{\partial{#1}}{\partial \mathcal{N}_i}}}
\newc{\pD}[2]{\dspl{\frac{\partial{#1}}{\partial #2}}}
\newc{\dt}{\dspl{\frac{d}{dt}}}
\newc{\bdry}[1]{\mbox{$\partial #1$}}
\newc{\sgn}{\mbox{sign}}

\newc{\Hess}[1]{\frac{\partial^2 #1}{\pdh z_i \pdh z_j}}
\newc{\hess}[1]{\partial^2 #1/\pdh z_i \pdh z_j}

%%function spaces

\newcommand{\Coone}[1]{\mbox{$C^{1}_{0}(#1)$}}
\newcommand{\lspac}[2]{\mbox{$L^{#1}(#2)$}}
\newc{\hspac}[2]{\mbox{$C^{0,#1}(#2)$}}
\newc{\Hspac}[2]{\mbox{$C^{1,#1}(#2)$}}
\newc{\Hosp}{\mbox{$H^{1}_{0}$}}
\newcommand{\Lsp}[1]{\mbox{$L^{#1}(\Og)$}}
\newc{\hsp}{\Hosp(\Og)}

%%greeks

\newc{\ag}{\alpha}
\newc{\bg}{\beta}
\newc{\cg}{\gamma}\newc{\Cg}{\Gamma}
\newc{\dg}{\delta}\newc{\Dg}{\Delta}
\newc{\eg}{\varepsilon}
\newc{\zg}{\zeta}
\newc{\thg}{\theta}
\newc{\llg}{\lambda}\newc{\LLg}{\Lambda}
\newc{\kg}{\kappa}
\newc{\rg}{\rho}
\newc{\sg}{\sigma}\newc{\Sg}{\Sigma}
\newc{\tg}{\tau}
\newc{\fg}{\phi}\newc{\Fg}{\Phi}
\newc{\vfg}{\varphi}
\newc{\og}{\omega}\newc{\Og}{\Omega}
%\newc{\ng}{\eta}
\newc{\pdh}{\partial}

%%Integration and Sum

\newc{\ii}[1]{\int_{#1}}
\newc{\iidx}[2]{{\dspl\int_{#1}~#2~dx}}
\newc{\bii}[1]{{\dspl \ii{#1} }}
\newc{\biii}[2]{{\dspl \iii{#1}{#2} }}
\newc{\su}[2]{\sum_{#1}^{#2}}
\newc{\bsu}[2]{{\dspl \su{#1}{#2} }}

\newcommand{\iiomdx}[1]{{\dspl\int_{\Og}~ #1 ~dx}}
\newc{\biiom}[1]{{\dspl\int_{\bdrom}~ #1 ~d\sg}}
\newc{\io}[1]{{\dspl\int_{\Og}~ #1 ~dx}}
\newc{\bio}[1]{{\dspl\int_{\bdrom}~ #1 ~d\sg}}
\newc{\bsir}{\bsu{i=1}{r}}
\newc{\bsim}{\bsu{i=1}{m}}

\newc{\iibr}[2]{\iidx{\bprw{#1}}{#2}}
\newc{\Intbr}[1]{\iibr{R}{#1}}
\newc{\intbr}[1]{\iibr{\rg}{#1}}
\newc{\intt}[3]{\int_{#1}^{#2}\int_\Og~#3~dxdt}
%%\newc{\itQ}[2]{\dspl{\int\hspace{-2.5mm}\int_{#1}~#2~dxdt}}
%%\newc{\mitQ}[2]{\dspl{\rule[1mm]{4mm}{.3mm}\hspace{-5.3mm}\int\hspace{-2.5mm}\int_{#1}~#2~dxdt}}

\newc{\itQ}[2]{\dspl{\int\hspace{-2.5mm}\int_{#1}~#2~dz}}
\newc{\mitQ}[2]{\dspl{\rule[1mm]{4mm}{.3mm}\hspace{-5.3mm}\int\hspace{-2.5mm}\int_{#1}~#2~dz}}
\newc{\mitQQ}[3]{\dspl{\rule[1mm]{4mm}{.3mm}\hspace{-5.3mm}\int\hspace{-2.5mm}\int_{#1}~#2~#3}}

\newc{\mitS}{\dspl{\rule[1mm]{3mm}{.3mm}\hspace{-4mm}\int}}
\newc{\mitx}[2]{\dspl{\rule[1mm]{3mm}{.3mm}\hspace{-4mm}\int_{#1}~#2~dx}}

\newc{\mitQq}[2]{\dspl{\rule[1mm]{4mm}{.3mm}\hspace{-5.3mm}\int\hspace{-2.5mm}\int_{#1}~#2~d\bar{z}}}
\newc{\itQq}[2]{\dspl{\int\hspace{-2.5mm}\int_{#1}~#2~d\bar{z}}}

\newc{\pder}[2]{\dspl{\frac{\partial #1}{\partial #2}}}

%%variables

\newc{\ui}{u_{i}}
\newcommand{\upl}{u^{+}}
\newcommand{\umn}{u^{-}}
\newcommand{\un}{\{ u_{n}\}}

\newcommand{\uo}{u_{0}}
\newc{\voi}{v_{i}^{0}}
\newc{\uoi}{u_{i}^{0}}
\newc{\vu}{\vec{u}}

\newc{\xo}{x_{0}}
\newc{\Br}{B_{R}}
\newc{\Bro}{\Br (\xo)}
\newc{\bdrom}{\bdry{\Og}}
\newc{\ogr}[1]{\Og_{#1}}
\newc{\Bxo}{B_{x_0}}

%%Additional macros
\newc{\inP}[2]{\|#1(\bullet,t)\|_#2\in\cP}
\newc{\cO}{{\mathcal O}}
\newc{\inO}[2]{\|#1(\bullet,t)\|_#2\in\cO}

\newc{\newl}{\\ &&}

%\newc{\bilhom}{\mbox{Bil}(\mbox{Hom}(\RR^n,\RR^N))}
\newc{\bilhom}{\mbox{Bil}(\mbox{Hom}(\RR^{nm},\RR^{nm}))}
\newc{\VV}[1]{{V(Q_{#1})}}

\newc{\ccA}{{\mathcal A}}
\newc{\ccB}{{\mathcal B}}
\newc{\ccC}{{\mathcal C}}
\newc{\ccD}{{\mathcal D}}
\newc{\ccE}{{\mathcal E}}
\newc{\ccH}{\mathcal{H}}
\newc{\ccF}{\mathcal{F}}
\newc{\ccI}{{\mathcal I}}
\newc{\ccJ}{{\mathcal J}}
\newc{\ccP}{{\mathcal P}}
\newc{\ccQ}{{\mathcal Q}}
\newc{\ccR}{{\mathcal R}}
\newc{\ccS}{{\mathcal S}}
\newc{\ccT}{{\mathcal T}}
\newc{\ccX}{{\mathcal X}}
\newc{\ccY}{{\mathcal Y}}
\newc{\ccZ}{{\mathcal Z}}

\newc{\bb}[1]{{\mathbf #1}}
\newc{\bbA}{{\mathbf A}}
%%Bracket
\newc{\myprod}[1]{\langle #1 \rangle}
\newc{\mypar}[1]{\left( #1 \right)}

%%Norms

\newc{\lspn}[2]{\mbox{$\| #1\|_{\Lsp{#2}}$}}
\newc{\Lpn}[2]{\mbox{$\| #1\|_{#2}$}}
\newc{\Hn}[1]{\mbox{$\| #1\|_{H^1(\Og)}$}}

%%Misc

\newc{\cyl}[1]{\og\times \{#1\}}
\newc{\cyll}{\og\times[0,1]}
\newc{\vx}[1]{v\cdot #1}
\newc{\vtx}[1]{v(t,x)\cdot #1}
\newc{\vn}{\vx{n}}

\newcommand{\RR}{{\rm I\kern -1.6pt{\rm R}}}

%\numberwithin{equation}{section}
%\newcommand{\eproof}{\vrule height5pt width3pt depth0pt}
%\newtheorem{exam}[lem]{Example}

\newenvironment{proof}{\noindent\textbf{Proof.}\ }
{\nopagebreak\hbox{ }\hfill$\Box$\bigskip}

%%End macros

\vspace*{-.8in}
\begin{center} {\LARGE\em Regularity for Fully Nonlinear P-Laplacian Parabolic Systems: the Degenerate Case}

 \end{center}

\vspace{.1in}

\begin{center}

{\sc Dung Le}{\footnote {Department of Mathematics, University of Texas at San
Antonio, One UTSA Circle, San Antonio, TX 78249. {\tt Email: dle@math.utsa.edu}
} }

\end{center}

\begin{abstract}
This paper studies H\"older
regularity property of bounded weak solutions to a class of strongly
coupled degenerate parabolic systems.
\end{abstract}

\vspace{.2in}

\section{Introduction}

We study the H\"older regularity of bounded weak solutions
of nonlinear $p$-Laplacian parabolic systems of the form

\beqno{e1} u_t = \Div(A(u,Du)) + F(u,Du), \eeq

\nid in a domain $Q=\Og\times(0,T)\subset \RR^{n+1}$, with $\Og$
being an open subset of $\RR^n$, $n\ge1$. The {\em vector valued}
function $u, f$ take values in $\RR^m$, $m\ge1$. $Du$ denotes the
spatial derivative of $u$. Here, $A(u,\zeta)$ is a nonlinear map from $\RR^m\times \RR^{nm}$ into $\RR^{nm}$.

A weak solution $u$ to \mref{e1} is a function $u\in
W^{1,0}_2(Q,\RR^m)$ such that $$\itQ{Q}{[\myprod{-u,\fg_t} + \myprod{A(u,Du),D\fg}]}
= \itQ{Q}{\myprod{F(u,Du),\fg}} $$ for all $\fg\in C_0^1(Q,\RR^m)$. Here,
we write $dz = dxdt$.

The evolution $p$-Laplacian scalar equation has been one of most widely studied nonlinear degenerate parabolic equations. The particular feature of \mref{e1} is its gradient-dependent diffusivity. Such systems, and their stationary counterparts, appear in different models in non-Newtonian fluids, turbulent flows in porous media, certain diffusion or heat transfer processes, and recently in image processing. 

A large body of literature on $p$-Laplacian systems has been devoted to the following system $$u_t = \Div(|Du|^{p-2}Du) + F(u,Du)$$ which is a special case of \mref{e1} where $A(u,Du)=|Du|^{p-2}Du$ does not depend explicitly on $u$. In this case, the regularity theory of bounded weak solutions was then almost settled and masterfully presented in the text book \cite{Dib} (see also \cite{Giusti} for the stationary counterpart). The techniques and results also hold for systems where $A$ depends smoothly on $x,t$. In fact, under suitable assumptions on $F$, we now know that bounded weak solutions to the above systems has H\"older continuous spatial derivatives. The theory was then based on a far-reaching combination of generalized DiGiorgi and Moser's methods for scalar equations.

However, this method breaks down in dealing with systems \mref{e1} allowing more general structural conditions and with the diffusivity $A$ depending explicitly on the unkown $u$. First of all, the dispersion of the eigenvalues of the derivative of $A$ with respect to the second variable $Du$ will prevent the Moser type iteration techniques in \cite[Chapter IX]{Dib} from being applicable in order to show that $|Du|$ is locally bounded, a starting and crucial point in defining the scaled cylinders in the next steps. Secondly, the presence of $u$ in $A$ will create extra terms when one differentiates the system in order to obtain a new system satisfied by $Du$. These extra terms may not be well defined if $u$ is not yet known to be H\"older continuous.

In this work, we choose a different approach. We will establish the H\"older continuity of $u$ by using a homotopy argument. We assume that the system \mref{e1} can be imbedded in a family of systems and at least one of which has the property that its bounded weak solutions are H\"older continuous and satisfy a {\em scaling decay estimate}. Under suitable assumptions, we show that this property will be carried onto bounded weak solutions of the considered system. This type of decay estimates with scaling was also used in \cite{Dib} using the local supremum norm of $|Du|$. In our case, since $|Du|$ is not locally bounded and the best we can say is that $|Du|^q$ is locally integrable for some $q>p$. The scaled cylinders in this work  must then be defined differently. We will use the average mean of $|Du|^p$ instead of its unavailable supremum norm.

Thus, we will consider a family of parabolic systems parameterized
by $\nu\in[0,1]$ \beqno{e1zz}u_t = \Div(A(\nu,u,Du)) + F(\nu,u,Du),
\mbox{ in $Q=\Og\times(0,T)\subset \RR^{n+1}$.}\eeq Assuming
\mref{e1} is the above system when $\nu=1$ and H\"older continuity
results are known for the system when $\nu=0$. Inspired by \cite[Proposition 3.1]{Dib}, we introduce the so call {\em scaling decay property} D) which is H\"older continuity for weak solutions to scalar degenerate equations. We then consider a subset $\ccI$
of parameters in $[0,1]$ where bounded weak solutions of the above system
satisfy this property. The main goal is then to prove that $\ccI$ is
both open and closed in $[0,1]$ so that $\ccI=[0,1]$ and the desired
H\"older continuity for solutions to \mref{e1} is obtained. 
Our first two main results concerning the open and closed properties
of $\ccI$ will be presented under two sets of conditions as they will
be established by using different tools, and they may be
independently of interest in other applications.

The main vehicle in the proof of $\ccI$ being open is the
\refprop{heat-aab}, which is the $p$-Laplacian
version of the nonlinear heat approximation result in \cite{ld5}.
Basically, it asserts that if a vector valued function $u$ {\em
almost} and weakly solves a system like \mref{e1zz}, with $\nu\in
I$, then it can be approximated in certain controllable way by a
solution $v$ of the system. By this, property D) of $v$
can be carried over to $u$.  The proof of this $p$-Laplacian
approximation version is not a simple extension of the result in \cite{ld5} as our systems are degenerate (or singular) and many more technical tools. Among them is a measure theoretic result \reflemm{weakconv} in \refsec{prem} establishing uniform continuity of the integrals of the derivatives of approximated solutions. As a consequence of this, in \reflemm{duLq}, we also present a result on higher integrability of the derivatives of "almost" weak solutions to a $p$-Laplacian system. Similar results for weak solutions to $p$-Laplacian systems were first reported in \cite{ki}.

On the other hand, the above argument is local by nature and cannot
be used to prove that $I$ is closed as it lacks certain uniform
estimates in order to show that limits of a sequence of regular
solutions are also regular. To this end, we will use a different
approach deriving uniform and global estimates for the integrals of
spatial derivatives of regular solutions with uniform bounded norms.

In this paper, our main results only concern the degenerate case, i.e. $p>2$. The singular case, $p<2$, can be dealt with in a similar way but much more subtle and will be reported in a forthcoming work. However, most of our main tools work for both cases and we report them here in \refsec{prem} and \refsec{app-sec} for future references. We will specifically state the range of $p$ for which our results hold.

The paper is organized as follows. In \refsec{mainsec}, we
will introduce notations and discuss in details our hypotheses and
main theorems. \refsec{prem} collects technical lemmas.
\refsec{app-sec} presents our main vehicles - the $p$-Laplacian nonlinear
approximation results. The proof that $I$ is open will be given in
\refsec{decaysec}. Finally, \refsec{iclosed} details the proof of $I$
being closed  and concludes our
paper.

\section{Main results}\eqnoset \label{mainsec}

Throughout this paper, for some $z_0=(x_0,t_0)\in \RR^{n+1}$ and
$R,\rg>0$, $Q_{R,\rg}(z_0)$ denotes the parabolic cylinder centered
at $z_0$ with radius $R,\rg$. That is,
$Q_{R,\rg}(z_0)=B_R(x_0)\times[t_0-\rg,t_0]$. We also abbreviate by
$S_{R,\rg}(z_0)=B_R(x_0)\times\{t_0-\rg\}\cup
\partial B_R(x_0)\times[t_0-\rg,t_0]$, the parabolic boundary of
$Q_R(z_0)$. If the center $z_0$ was understood, we would simply
write $Q_{R,\rg},S_{R,\rg}$ for $Q_{R,\rg}(z_0)$ and
$S_{R,\rg}(z_0)$ respectively.

For a given cylinder $Q=B\times[a,b]$ and $p>1$, we consider the
space $V_p(Q)=V_p(Q,\RR^m)$ of vector valued functions $u:Q\to\RR^m$
with norm $\|\cdot\|_{V_p(Q)}$ defined by
$$\|u\|_{V_p(Q)} = \sup_{t\in[a,b]}\|u(\cdot,t)\|_{L^2(B)} + \|Du\|_{L^p(Q)}.$$ By
$V_p^0(Q)$ we denote the closure of $C_0^1(Q)$ in $V_p(Q)$ with
respect to the above norm.

Let $z_0=(x_0,t_0)$ and $Q_{R,\rg}(z_0)$ be any parabolic cylinder in $R^{n+1}$, the
following scaled norm will also be used $$\|u\|_{V(Q_{R,\rg})} =
\sup_{t\in [t_0-\rg,t_0]}\left(R^{-n}\iidx{B_R(x_0)}{u(x,t)}\right)^\frac{1}{2} + \left(R^{p}R^{-n}\rg^{-1}\itQ{Q_{R,\rg}}{|Du|^p}
\right)^\frac{1}{p}.$$ Obviously, this norm is invariant via dilations.

For any integrable function $u:Q\to\RR^m$ and any measurable subset
$A$ of $Q$, we write $$ u_A = \frac{1}{|A|}\int_A u(z)dz =
\mitQ{A}{u(z)}.$$ If $A$ is a cylinder $Q_{R,\rg}=Q_{R,\rg}(z_0)$
and there is no possibility of ambiguity, we simply write
$u_{R,\rg}=u_{Q_{R,\rg}}$. Furthermore, if $\rg$ is defined in term of $R$ and the relation between
$R,\rg$ is clear we also abbreviate  $u_{R,\rg}$ by $u_R$ for the
sake of simplicity.

As our results are local in nature, without loss of generality, we
will simply consider $Q$ being the unit parabolic cylinder
$B_1(0)\times[-1,0]$ throughout this paper. We then consider a
family of systems \beqno{fame1a} u_t = \Div(A(\nu,u,Du)) +
F(\nu,u,Du), \mbox{ in $Q=B_1\times(-1,0)\subset \RR^{n+1}$ and
$\nu\in[0,1]$.}\eeq

By a bounded weak solution $u$ to this system we mean a bounded
vector valued function $u$ satisfying
$$ \itQ{Q}{[-u\fg_t+\myprod{A(\nu,u,Du),D\fg}]} -\itQ{Q}{F(\nu,u,Du)\fg} =0, \quad \forall \fg\in
C^1_0(Q,\RR^m).$$

For simplicity, we will mainly consider the case $F\equiv0$ in our discussion. The presence of $F$ can be treated with minor modifications and we will briefly discuss this case at the end of this section. 

We will always consider matrices $A(\nu,u,\zeta)$ satisfying the
following ellipticity condition

\bdes \item[E)] There are positive constants $\llg,\LLg$ such that
for all $\zeta,\eta\in\RR^{nm},u\in \RR^m$ \beqno{aaa}
\myprod{A(\nu,u,\zeta),\zeta}\ge\llg|\zeta|^p, \quad
|\myprod{A(\nu,u,\zeta),\eta}|\le\LLg|\zeta|^{p-1}|\eta| .\eeq \edes

In the study of the H\"older regularity of  a weak solution $u$, it
is now well known that (see \cite{GiaS}) one needs to establish a
{\em mean oscillation decay estimate}: For some $\tau,\ag\in(0,1)$
and any $Q_R=Q_{R,\rg}(x,t)\subset Q$, there are positive constant $\tau_0, C(\tau_0)$ such
that \beqno{decay} \mitQ{Q_{\tau R,\tau\rg}}{|u-u_{\tau
R,\tau\rg}|^2} \le C(\tau_0)\tau^\ag\mitQ{Q_{R,\rg}}{|u-u_{R,\rg}|^2} \quad \forall \tau\in(0,\tau_0).\eeq

However, due to the degeneracy/singularity of the diffusion matrix
$A$, such decay estimates do not hold in general for cylinders whose
space-time configuration are uniform for all solutions and
depend only on the parameters defining the systems. Roughly speaking, we
can only obtain here \mref{decay} when $R,\rg$ are linked via an
intrinsic scaling determined by the solution $u$ itself. Yet this
decay property still gives the desired H\"older continuity. The idea of using scaled cylinders was known in literature (see \cite{Dib}) where scalar equations were studied so that Harnack type inequalities could be established and gave the H\"older continuity. Here, such techniques are no longer available and we have to scale the cylinders by using the average oscillations of weak solutions. Being inspired by \cite[Proposition 3.1]{Dib}, we introduce the following decay property. Note that the mean oscillation of vector valued solutions is used here in place of the essential oscillation for scalar solutions in \cite{Dib}. 

We say that a bounded vector valued function $v:Q\to\RR^m$
satisfies a {\em scaling decay property} if

\bdes \item[D)] Let $M=\sup_Q |v|$. For any $R_0>0, 
\eta\in(0,1)$, and $(x_0,t_0)\in Q$ there are
 positive numbers $A,K,L,\ag_0,\og_0$ depending on
$M,\eta$ (with $K,A$ sufficiently large) such that we can define the
following sequences\beqno{scalezz} R_k=\frac{R_0}{K^k}, \quad
\og_{k+1}=\max\{\eta\og_k, LR_{n}^{\ag_0}\},\quad
S_k=\frac{\og_k}{A}, Q_k=B_{R_k}(x_0)\times
[t_0-S_k^{2-p}R_k^p,t_0].\eeq such that if  $$\og_0^p \ge \mitQ{Q_0}{|v-(v)_0|^p},\quad
(v)_0=\mitQ{Q_0}{v} $$ then for any integer $k=1,\ldots$
\beqno{vazz} \og_k^p \ge \mitQ{Q_k}{|v-(v)_k|^p},\quad
(v)_k=\mitQ{Q_k}{v}.\eeq \edes

Our first main result shows that if \mref{fame1a}, for some $\nu$,
is a "nice" system in the sense that its bounded weak solutions
satisfying the decay property D), then "near by" systems are also
nice.
To be more precise, let us describe the this result in details here.
We first suppose that the family of systems \mref{fame1a} contains
at least a "nice" one.

\bdes

\item[I)] There is a nonempty set $\ccI\subset[0,1]$ such that for
any $\nu$ in $I$ the decay property D) holds for any bounded weak
solution to \mref{fame1a}. The same assumption applies to the systems  with $u$ being replaced by any constant vector $C$, i.e.  $$ \itQ{Q}{[-u\fg_t+\myprod{A(\nu,C,Du),D\fg}]} -\itQ{Q}{F(\nu,C,Du)\fg} =0, \quad \forall \fg\in
C^1_0(Q,\RR^m)$$
\edes

We then consider the following structural assumptions on the
matrices $A(\nu,u,Du)$.

\bdes

\item[O.1)] (Uniform ellipticity) For any $\nu\in[0,1]$, $A(\nu,u,\zeta)$
satisfies the ellipticity condition E) for some positive constants
$\llg,\LLg$. In addition, $\pder{A}{\zeta}(\nu,u,\zeta)$ exists and there is $c_0>0$ such that $$
\myprod{\pder{A}{\zeta}(\nu,u,\zeta),\zeta}\ge c_0|\zeta|^p$$ for any bounded weak solution $u$ of 
\mref{fame1a}.

\item[O.2)]  (Monotonicity) For any $w\in \RR^m$ and $U,V\in\RR^{nm}$, there
holds \beqno{AAAmono} \myprod{A(w,U)-A(w,V),U-V} \ge
\llg_0(U,V)|U-V|^2, \eeq where
$$\llg_0(U,V)=c_0\left\{\barr{ll}\min\{|U|^{p-2},|V|^{p-2}\}&
U\ne0\mbox{ or } V\ne0\\0&\mbox{otherwise}\earr\right.$$ for some
positive constant $c_0$.

\item[O.3)] (Continuity) $A(\nu,u,\zeta)$  is
H\"older in $\nu\in[0,1]$ and Lipschitz in $u$. That is,
\beqno{Aconvz} |A(\nu,u,\zeta)-A(\mu,u,\zeta)| \le
C|\nu-\mu|^\thg|\zeta|^{p-1}, \eeq \beqno{ak2a}
|A(\nu,u,\zeta)-A(\nu,v,\zeta)|\le C|u-v||\zeta|^{p-1}\eeq for some
$C,\thg>0$ and any $\nu,\mu\in[0,1]$, $u,v\in \RR^m,\,
\zeta\in\RR^{nm}$.

\item[O.4)] (Existence)  For $\nu\in \ccI$ and any cylinder
$Q'\subset Q$, and any vector valued function $g\in C^1(Q')$, the
system \beqno{exist}\left\{\barr{l}
\itQ{Q_R}{[-u\fg_t+\myprod{A(\nu,u,Du),D\fg}]} =0, \quad \forall
\fg\in C^1_0(Q',\RR^m),
\\ \mbox{ $u=g$ on $S_{Q'}$} \earr\right.\eeq has a bounded weak
solution $u$ (and thus satisfies the property D)).
\item[O.5)] (Uniform maximum principle) For any $\llg\in[0,1]$ and $Q'\subset Q$, if $g\in C^1(Q')$ then
there is a constant $C(\|g\|_{L^\infty(Q')})$ such that any weak
solution $u$ to the system \mref{exist} is bounded and satisfies the
estimate $\|u\|_{L^\infty(Q')} \le C(\|g\|_{L^\infty(Q')})$.

\edes

\brem{h4rem} We should remark that  if the argument $u$ in $A(\nu,u,Du)$ is replaced by a constant vector and the data $g$ is bounded then the existence of a weak solution in O.4) is granted by classical approximation methods as in \cite{LSU} or \cite{Lions}. Otherwise, the existence condition O.4) can be satisfied by using Gal\"erkin's method if $g$ is sufficiently regular and the solution $u$ is known a-priori to be H\"older continuous. Since $\ccI$ is the set of parameters for which the systems has bounded weak solutions being H\"older continuous, O.4)  is justified.\erem

Our first main result then asserts that the set $\ccI$, where the
property D) holds, is open.

\btheo{regmain} Suppose that I) and O.1)-O.5) hold and $p>2$. Then $\ccI$ is
open in the usual topology of $[0,1]$. Moreover, bounded weak solutions to the system \mref{fame1a} with $\nu\in\ccI$ are H\"older continuous. \etheo

The above theorem relies on a nontrivial generalization of the so called nonlinear heat approximation lemma which was introduced in our earlier work \cite{ld5} concerning nondegenerate systems.

Next, we will give conditions for the set $\ccI$ to be closed in
$[0,1]$. To this end, we take a sequence $\{\nu_k\}$ in $\ccI$ such
that $\nu_k\to\mu$ and we will show that $\mu\in \ccI$. We first
require that any bounded weak solution $u$ to \mref{fame1a}, with
$\nu=\mu$, can be weakly approximated by "nice" solutions.

\bdes \item[II)] For each $\nu\in \ccI$, the system \mref{fame1a}
satisfies the existence condition O.4) and maximum principle O.6).
Moreover, if ${\nu_k}\subset \ccI$ and $\nu_k\to\mu$ then for any
bounded weak solution $u$ to \mref{fame1a} with $\nu=\mu$ there is a
sequence $\{v_k\}$ of H\"older continuous solutions to \mref{fame1a}, with
$\nu=\nu_k$, such that $Dv_k$ converges weakly to $Du$ in $L^p(Q)$
and the $L^\infty$ norms of $v_k$ are bounded uniformly in terms of
that of $u$.
 \edes

Although it will be shown in \refsec{app-sec} that the above
assumption holds under O.1)-O.5) via nonlinear heat approximation, we state II) here for our next
result, which can be of interest in itself, so that it is
independent of \reftheo{regmain}. Apparently, II) could also be
verified by other means via weaker assumptions than O.1)-O.5).

We then consider $\nu\in \ccI$ and a $C^1$ solution $v$ to
\beqno{s1dvz} \left\{\barr{ll}v_t = \Div(A(\nu,v,Dv)) & \mbox{ in
$Q$,}\\ v=g & \mbox{ on $S$.}\earr\right. \eeq

The boundary condition $g$ is assumed to be smooth. By II), the
above system satisfies the maximum principle and we can define
$$M_{\nu,v}=\sup_{Q_\frac34}|v|.$$

For any bounded weak solution $v$ to \mref{s1dvz}, we then impose the following assumptions on the structure of the system.

\bdes
\item[M.1)] The matrix
$(A^{ij}_{kl})=\pder{A}{\zeta}(\nu,v,\zeta)$ is elliptic with
 the ellipticity constants  $\llg_{\nu,v},\LLg_{\nu,v}$.
That is \beqno{A-ella}\sum_{i,j=1}^m\sum_{k,l=1}^n
A^{ij}_{kl}\eta^i_k\eta^j_l \ge \llg_{\nu,v}|\eta|^2, \quad
\sum_{i,k}(\sum_{j,l}A^{ij}_{kl}\eta^j_l)^2\le
\LLg_{\nu,v}^2|\eta|^2 \quad \forall \eta\in\RR^{nm}.\eeq Moreover, for some positive constants
$\llg_\nu,\LLg_\nu$ we have
\beqno{lambda-z}\llg_{\nu,v}\ge \llg_\nu|\zeta|^{p-2}, \quad \LLg_{\nu,v}\le
\LLg_\nu|\zeta|^{p-2}.\eeq If $n>2$, we also assume that
\beqno{lambda-conda}
\sup\{\frac{\LLg_{\nu,v}}{\llg_{\nu,v}}\,:\,\mbox{ $v$ is a bounded weak
solution to \mref{s1dvz}}\}<\frac{n}{n-2}.\eeq
\item[M.2)]  For every $\nu$ in $I$ and any bounded weak solution $v$ to \mref{s1dvz},
there exists a positive constant $a_{\nu,v}$ such that \beqno{MMa}
|\frac{\partial A} {\partial v} (\nu,v,\xi)| \le
a_{\nu,v}|\xi|^{p-1} \mbox{ with } 2a_{\nu,v}
M_{\nu,v}(p+n-1)<\sg_0\widehat{\llg}_{\nu}, \eeq where $\sg_0$ is a
fixed number in $(0,1)$ and
$$\widehat{\llg}_{\nu} = (1-\dg^2)\llg_{\nu} \mbox{ and }
\dg=\frac{n-2}{n}\sup\{\frac{\LLg_{\nu,v}}{\llg_{\nu,v}}\,:\,\mbox{
$v$ is a bounded weak solution to \mref{s1dvz}}\}.$$

\edes

Note that $\widehat{\llg}_{\nu,v}>0$ due to \mref{lambda-conda}.
Meanwhile, \mref{lambda-conda} requires that the principal
eigenvalues $\LLg_{\nu,v},\llg_{\nu,v}$ of $\pder{A}{\zeta}$ are not
too far apart (when $n>2$). We should also remark that the constants $\llg_\nu,\LLg_\nu$ could be allowed to depend on $v$ as long as there were fixed positive numbers $c_1,c_2$ such that  quotient $c_1\le \llg_\nu/\LLg_\nu\le c_2$. We assume however that they are constants for the sake of simplicity. 

Our
next main result reads

\btheo{dvp1} Assume  the conditions II) and M.1)- M.2) and that $p>2$. The set $I$
is closed in [0,1]. \etheo

Combining with the results of the previous theorem, as O.1)-O.5) are
sufficient for II), we then have

\btheo{dvp1z} Assume  I), O.1)-O.5) and M.1)- M.2) and $p>2$ then $I=[0,1]$. Thus, bounded weak solutions are H\"older continuous.
\etheo

Finally, we remark that our proof continues to hold for systems like \beqno{fame1ab} u_t = \Div(A(\nu,u,Du)) +
F(\nu,u,Du),\eeq if the nonlinearity $F$ satisfies
$$|\pder{F}{v}(\nu,v,\xi)| \le C + C|v|^l +\eg_0|\xi|^p, \quad
|\pder{F}{\xi}(\nu,v,\xi)| \le C + C|v|^l +\eg_0|\xi|^{p-1}$$ for some $l>0$ and
sufficiently small $\eg_0>0$. The presence of $F$ would cause extra terms in our arguments but they can be easily treated by invoking the H\"older and Young inequalities.

\section{Technical lemmas}\label{prem}\eqnoset
In this section, we present various estimates on a vector valued
function $u$ weakly satisfying certain differential inequality. Although our main results in this work concern the degenerate situation when $p>2$ and the singular case ($p<2$) will be treated in future works,
several results in this section hold for $p>1$ and we will specify
the range of $p$ in each statement for future reference.

Throughout this paper, the constants $C,C_1,\ldots$ can change line by line but they are all universal constants as they depend only on the initially fixed parameters (such as $n,m$). For any two quantities $A,B$, we write $A\sim B$ if there are universal positive constants $C_1,C_2$ such that $C_1A\le B\le C_2A$.

First of all, we recall the following Sobolev inequalities in a ball of $\RR^n$
and the dependence of the constants on the size of the ball. Let
$\fg$ be a function defined on $B_1$. By scaling, with
$x=\frac1R\bar{x}$ and $\fg(x)=\bar{\fg}(\bar{x})$, we have  the followings facts on the $k$th derivatives $D^{(k)}$
$$\|D_x^{(k)}\fg\|^p_{L^{p}(B_1)} \sim
CR^{pk-n}\|D_{\bar{x}}^{(k)}\bar{\fg}\|^p_{L^{p}(B_R)} \mbox{ and
}\|D_x\fg\|^p_{L^{q}(B_1)} \sim
CR^{(1-\frac{n}{q})p}\|D_{\bar{x}}\bar{\fg}\|^p_{L^{q}(B_R)}.$$ In
particular, for $q=p_*=pn/(n+p)$, we easily see that
$\|D_x\fg\|^p_{L^{q}(B_1)} \sim
CR^{-n}\|D_{\bar{x}}\bar{\fg}\|^p_{L^{q}(B_R)}$. Using this in the
Sobolev-Poincar\'e inequality on $B_1$, we get
$$\|\fg\|^p_{L^{p}(B_1)}\le C(n)\|D_x\fg\|^p_{L^{q}(B_1)}\Rightarrow
\|\bar{\fg}\|^p_{L^{p}(B_R)}\le
C(n)R^{n+(1-\frac{n}{q})p}\|D_{\bar{x}}\bar{\fg}\|^p_{L^{q}(B_R)}$$
for all $p$ such that $q\le p \le q^*=qn(n-q)$. 

On the other hand,
$$\|\fg\|_{W^{-k,p'}(B_1)} = \sup_\psi\frac{1}{\|\psi\|_{W_0^{k,p}(B_1)}}\int_{B_1}\fg\psi
dx\le\sup_\psi\frac{C}{R^\frac{pk-n}{p}\|\bar{\psi}\|_{W_0^{k,p}(B_R)}}\frac{1}{R^n}\int_{B_R}\bar{\fg}\bar{\psi}
d\bar{x}.$$

This also implies $\|\fg\|^p_{W^{-k,p'}(B_1)} \sim
CR^{-pk+n-pn}\|\bar{\fg}\|^p_{W^{-k,p'}(B_R)}$.

We begin with the following lemma.

\blemm{pineq} Let $Q_{R,\rg}=B_R\times(-\rg,0)$ and
$u:Q_{R,\rg}\to\RR^m$ be in $V_p(Q_{R,\rg})$ for some $p>1$. Assume that
there is a function $G\in L^1(Q_{R,\rg})$ such that
\beqno{p1} \left|\itQ{Q_{R,\rg}}{u\fg_t}\right| \le
C\itQ{Q_{R,\rg}}{|G||D\fg|} \quad \forall \fg\in
C^1_0(Q_{R,\rg}).\eeq
 Let $q$ be such
that $p_*=np/(n+p)\le q < p$. Then, for any $\eg>0$ there exist
positive constants $C, C(\eg)$ such that
\beqno{pp}\barr{lll}\mitQ{Q_{R,\rg}}{|u-u_R|^p} &\le& \eg
R^p\mitQ{Q_{R,\rg}}{|Du|^p} +
(C+\eg)R^{p}\left(\mitQ{Q_{R,\rg}}{|Du|^q}\right)^\frac{p}{q}\\
&&+C(\eg)(\frac{\rg}{R})^{p}\left(\mitQ{Q_{R,\rg}}{|G|}\right)^p.\earr\eeq
 \elemm

\bproof For $s,r\in (-\rg,0)$ and $\eg>0$, take $\fg =
\psi(x)\eta(t)$ where $\psi\in C^1_0(B_R)$ and $\eta\equiv 1$ in
$(s,r)$, $\eta$ is linear in $(s-\eg,s)$ and $(r,r+\eg)$, $\eta$ is
zero elsewhere. By the Young inequality, the above estimate becomes
$$\left|\dspl{\int_{B_R}\left(\frac1\eg\int_{s-\eg}^sudt
-\frac1\eg\int_r^{r+\eg}udt\right)\psi dx}\right|\le
C\left(\itQ{Q_{R,\rg}}{|G\eta|}\right)\|\psi\|_{C_0^1(B_R)}.$$

Let $l>(n+p)/p$ and $m=l-1>n/p$. By scaling, with $x=\frac1R\bar{x}$
and $\fg(x)=\bar{\fg}(\bar{x})$, we have for $\fg\in W_0^{m,p}(B_1)$
that
$$\|\fg\|_{L^{\infty}(B_1)}\le
C(n)\|D^{(m)}_x\fg\|_{L^{p}(B_1)}\Rightarrow
\|\bar{\fg}\|_{L^{\infty}(B_R)}\le
CR^{\frac{pm-n}{p}}\|D_{\bar{x}}^{(m)}\bar{\fg}\|_{L^{p}(B_R)},$$

Using this for $\bar{\fg}=D\psi$, we obtain $\|\psi\|_{C^1_0(B_R)}
\le CR^\frac{pl-p-n}{p}\|\psi\|_{W_0^{l,p}(B_R)}$. Hence, letting
$\eg\to0$, we obtain  the following estimate for $s,r\in (-\rg,0)$ $$\left|\int_{B_R}
(u(\cdot,s) -u(\cdot,r))\psi(x)dx\right| \le
CR^\frac{pl-p-n}{p}\left(\itQ{Q_{R,\rg}}{|G|}\right)\|\psi\|_{W_0^{l,p}(B_R)}
\quad \forall \psi\in W_0^{l,p}(B_R).$$ 
 Setting
$H(t,s,\cdot)=u(\cdot,s) -u(\cdot,t)$, we just proved that
\beqno{HHest}\|H(t,s,\cdot)\|_{W^{-l,p'}(B_R)}\le
CR^\frac{pl-p-n}{p}\left(\itQ{Q_{R,\rg}}{|G|}\right) \quad \forall s,t\in (-\rg,0).\eeq

Let $q$ be such that $p_*=np/(n+p)\le q < p$. We have $W^{1,q}(B_1)
\subset L^p(B_1)\subset W^{-l,p'}(B_1)$. Because $W^{1,q}(B_1)$ is
compactly imbedded in $L^p(B_1)$, a simple argument by contradiction
gives the following interpolation inequality
$$\|\fg\|_{L^p(B_1)} \le \eg\|\fg\|_{W^{1,q}(B_1)} +
C(\eg)\|\fg\|_{W^{-l,p'}(B_1)}, \quad \forall \fg\in W^{1,q}(B_1).$$

Using the norm $\|\fg\|_{W^{1,q}(B_1)} = \|D\fg\|_{L^{q}(B_1)} +
\|\fg\|_{L^{p}(B_1)}$ and choosing $\eg$ small, we have
$$\|\fg\|_{L^p(B_1)} \le \eg\|D\fg\|_{L^{q}(B_1)} +
C(\eg)\|\fg\|_{W^{-l,p'}(B_1)}, \quad \forall \fg\in W^{1,q}(B_1).$$

By a simple scaling argument, we derive from the above that
$$R^{-n}\|\fg\|^p_{L^p(B_R)} \le \eg R^{(1-\frac{n}{q})p}\|D\fg\|^p_{L^{q}(B_R)} +
C(\eg)R^{-lp+n-pn}\|\fg\|^p_{W^{-l,p'}(B_R)}, \quad \forall \fg\in
W^{1,q}(B_R).$$

Applying this to $H(t,s,\cdot)$ and using \mref{HHest}, we have
 \beqno{Hest}R^{-n}\|H(t,s,\cdot)\|^p_{L^p(B_R)} \le \eg
R^{(1-\frac{n}{q})p}\|DH(t,s,\cdot)\|^p_{L^{q}(B_R)} +
C(\eg)R^{-p-pn}\left(\itQ{Q_{R,\rg}}{|G|}\right)^p.\eeq

We now choose $s$ such that $$\int_{B_R}|Du(x,s)|^q\,dx \le
\frac{1}{\rg}\int_{-\rg}^0\int_{B_R}|Du(x,t)|^q\,dxdt.$$ Then
\beqno{schoice}\|Du(\cdot,s)\|^p_{L^q(B_R)} \le
\left(\frac{1}{\rg}\int_{-\rg}^0\int_{B_R}|Du(x,t)|^q\,dxdt\right)^\frac{p}{q}=\rg^{-p/q}\left(\itQ{Q_{R,\rg}}{|Du|^q}\right)^\frac{p}{q}.\eeq

On the other hand, by Sobolev-Poincar\'e's inequality, we also have
$$R^{-n}\|u(\cdot,s)-u(\cdot,s)_R\|^p_{L^p(B_R)} \le CR^{(1-\frac nq)p}\|Du(\cdot,s)\|^p_{L^q(B_R)} \le
CR^{p}\left(\mitQ{Q_{R,\rg}}{|Du|^q}\right)^\frac{p}{q}.$$

Obviously, $$\|u(\cdot,t)-u(\cdot,t)_R\|^p_{L^p(B_R)} \le
C\|u(\cdot,t)-u(\cdot,s)\|^p_{L^p(B_R)}+
\|u(\cdot,s)-u(\cdot,s)_R\|^p_{L^p(B_R)}$$ and $|DH(t,s,\cdot)|\le
|Du(\cdot,t)|+|Du(\cdot,s)|$. A simple use of H\"older's inequality
gives $\|Du(\cdot,t)\|^p_{L^q(B_R)} \le CR^{n(\frac pq-1)}
\|Du(\cdot,t)\|^p_{L^p(B_R)}$.

Together, when such $s$ is fixed, the above yields

$$\barr{lll}\frac{1}{R^n}\|u(\cdot,t)-u(\cdot,t)_R\|^p_{L^p(B_R)} &\le&
\eg R^{p-n}\|Du(\cdot,t)\|^p_{L^{p}(B_R)} + (C+\eg)CR^{p}\left(\mitQ{Q_{R,\rg}}{|Du|^q}\right)^\frac{p}{q}\\
&&+C(\eg)R^{-p-pn}\left(\itQ{Q_{R,\rg}}{|G|}\right)^p.\earr$$

Integrating the above over $t\in[-\rg,0]$ and dividing by $\rg$, we
get \mref{pp} and the proof of the lemma is then complete. \eproof

The following Poincar\'e type inequality is an immediate consequence
of the above lemma.

\blemm{reverse1} Assume as in \reflemm{pineq}. If $\rg=S^{2-p}R^2$
and $|G|\le |Du|^{p-1}$ then for any $\eg>0$ and $np/(n+p)\le q < p$ there exist
positive constants $C, C(\eg)$ such that
\beqno{lpkey1}\barr{lll}\mitQ{Q_{R,\rg}}{|u-u_R|^p} &\le& \eg
R^p\mitQ{Q_{R,\rg}}{|Du|^p} +
(C+\eg)R^{p}\left(\mitQ{Q_{R,\rg}}{|Du|^q}\right)^\frac{p}{q}\\
&&+C(\eg)(S^{2-p}R)^{p}\left(\mitQ{Q_{R,\rg}}{|Du|^{p-1}}\right)^{p}.\earr\eeq

\elemm

We now consider a weak solution $u\in V_p(Q_{1,1})$ to \beqno{plemeqn} u_t
= \Div(A(u,Du) \quad \mbox{in $Q_{1,1}$}. \eeq The matrix $A$ is assumed
to satisfy the following ellipticity conditions E) for some positive
constants $\llg,\LLg$ and $p>1$.

By testing \mref{plemeqn} with $|u-c|\fg^2$, with $\fg$ being a cutoff function
for $Q_{\frac12R,\frac12\rg}, Q_{R,\rg}$, we easily get the
following Caccioppoli type inequality

\blemm{CACCineq} Let $u$ satisfy \mref{plemeqn}. For any constant
vector $c\in\RR^m$ and any $Q_{R,\rg}\subset Q_{1,1}$ \beqno{rev-cacc}
\mitQ{Q_{\frac12R,\frac12\rg}}{|Du|^p} \le C
\frac{1}{R^p}\mitQ{Q_{R,\rg}}{|u-c|^p}+C\frac{1}{\rg}\mitQ{Q_{R,\rg}}{|u-c|^2}.\eeq
\elemm

A consequence of \reflemm{pineq} and the above is the following reverse H\"older
inequality.

\blemm{reverse2} Assume that $p>2$ and $u$ satisfies \mref{plemeqn}.
If $\rg=S^{2-p}R^2$ for some $S^p\sim \mitQ{Q_{R,\rg}}{|Du|^p}$ then we
have for any given positive $\eg$ and $np/(n+p)\le q < p$\beqno{revlp} \barr{lll}\mitQ{ Q_{\frac12R,\frac12\rg}}{|Du|^p}
&\le& \eg \mitQ{Q_{R,\rg}}{|Du|^p} +
(C+\eg)\left(\mitQ{Q_{R,\rg}}{|Du|^q}\right)^\frac{p}{q}\\
&&+C(\eg)\left(\mitQ{Q_{R,\rg}}{|Du|^{p-1}}\right)^\frac{p}{p-1}.\earr\eeq\elemm

\bproof For $\rg=S^{2-p}R^2$ we have \mref{lpkey1}. By H\"older
inequality we estimate the last term in \mref{lpkey1} as follows
$$(S^{2-p}R)^{p}\left(\mitQ{Q_{R,\rg}}{|Du|^{p-1}}\right)^{p}
\le R^p\left(\frac{\mitQ{Q_{R,\rg}}{|Du|^{p}}}{S^p}\right)^{(p-2)}
\left(\mitQ{Q_{R,\rg}}{|Du|^{p-1}}\right)^\frac{p}{p-1}.$$

If $p>2$ then the Young inequality ($1-(p-2)/p=2/p$) can apply to the last term in
\mref{rev-cacc} to yield $$\frac{1}{\rg}\mitQ{Q_{R,\rg}}{|u-u_R|^2}
= S^{p-2}R^{-2}\mitQ{Q_{R,\rg}}{|u-u_R|^2} \le \eg S^p +
C(\eg)\frac{1}{R^p}\mitQ{Q_{R,\rg}}{|u-u_R|^p}.$$

Thus, if $S^p\sim \mitQ{Q_R}{|Du|^p}$, we can combine the above
estimates with \mref{rev-cacc} and \mref{lpkey1} to obtain
\mref{revlp} if $p>2$. \eproof

The above result  also holds for $2n/(n+2)<p<2$ but we have to treat the
last term in \mref{rev-cacc} differently.

\blemm{rev-pless2} Assume that $2n/(n+2)<p<2$ and $u$ satisfies
\mref{plemeqn}. If $\rg=S^{2-p}R^2$ for some $S$
satisfying\beqno{S1} S^p \ge
c(n)\left( \mitQ{Q_{2R,2\rg}}{|Du|^p} +
\frac{1}{\rg}\mitQ{Q_{2R,2\rg}}{|u-u_{2R}|^2}\right)\eeq then we
have for any given positive $\eg$ and $q=2n/(n+2)$\beqno{revlpz} \barr{lll}\mitQ{ Q_{\frac12R,\frac12\rg}}{|Du|^p}
&\le& \eg \mitQ{Q_{R,\rg}}{|Du|^p} +
(C+\eg)\left(\mitQ{Q_{2R,2\rg}}{|Du|^q}\right)^\frac{p}{q}\\
&&+C(\eg)S^p.\earr\eeq\elemm

\bproof Let $\chi(x)$ be a smooth function in $x$ with compact
support in $B_R(x_0)$ such that $|D\chi|\le c/R$. We define
$$ u_{x_0,R}^\chi(t) = \int_{B_{R}(x_0)} u(x,t)dx /\int_{B_{R}(x_0)}
\chi(x)dx.$$ For any $s,t$ such that $-\rg<s<t<0$, we test \mref{plemeqn} with $(u_{x_0,R}^\chi(t)-u_{x_0,R}^\chi(s))\chi(x)\eta(\tau)$, $\eta$ is defined as in
\reflemm{pineq}, to get the following.

 $$|u_{x_0,R}^\chi(t)-u_{x_0,R}^\chi(s)|^2 \le
 C\frac{t-s}{R^{n+2}}\int_{Q_{R,\rg}}|G|^{2}dz
\le C\frac{t-s}{R^{n+2}}\int_{Q_{R,\rg}}|Du|^{2p-2}dz.$$

Thus, if $|t-s|\le\rg=S^{2-p}R^2$ then
\beqno{ww1}\frac{1}{\rg}|u_{x_0,R}^\chi(t)-u_{x_0,R}^\chi(s)|^2 \le
CS^{2-p}\mitQ{Q_{R,\rg}}{|Du|^{2p-2}}.\eeq

Obviously, for any $s\in(-\rg,0)$ $$\barr{lll}\frac{1}{\rg}\mitQ{Q_{R,\rg}}{|u-u_{R}|^2} &\le&\frac{1}{\rg}\mitQ{Q_{R}}{|u-u_{x_0,R}^\chi(s)|^2}\\&\le&
\frac{1}{\rg}\mitQ{Q_{R}}{|u-u_{x_0,R}^\chi(t)|^2} +
\frac{1}{\rg}\sup_{t\in(-\rg,0)}|u_{x_0,R}^\chi(t)-u_{x_0,R}^\chi(s)|^2.\earr$$ The last term can
be estimated as in \mref{ww1} while the first term on the right is bounded by$$ \itQ{Q_{R,\rg}}{|u-u_{x_0,R}^\chi(t)|^2}
\le \mypar{\sup_t\int_{B_{R}}|u-u_{x_0,R}^\chi(t)|^2dx}^{1-\frac
q2}\int_\rg^0\mypar{\int_{B_{R}}|u-u_{x_0,R}^\chi(t)|^2dx}^{\frac
q2}dt,$$ where  $q=2n/(n+2)$. The last factor can be bounded via Poincar\'e-Sobolev's inequality (in
the $x$ variable) by
$$\int_\rg^0\mypar{\int_{B_{R}}|u-u_{x_0,R}^\chi(t)|^2dx}^{\frac
q2}dt \le \itQ{Q_{R,\rg}}{|Du|^q}.$$ 

From the equation \mref{plemeqn} for $u$, we easily obtain

\beqno{uestpless2}\barr{lll}\lefteqn{\sup_{t\in(-\rg,0)}\int_{B_{R}}|u-u_{x_0,R}^\chi(t)|^2dx
\le \sup_{t\in(-\rg,0)}\int_{B_{R}}|u-u_{2R}|^2dx}\hspace{2cm}&&\\&\le&
\frac{1}{\rg}\itQ{Q_{2R,2\rg}}{|u-u_{2R}|^2} +
\frac{1}{R^p}\itQ{Q_{2R,2\rg}}{|u-u_{2R}|^p}.\earr\eeq

Applying H\"older inequality to the right hand side of the Poincar\'e inequality \mref{lpkey1}, we have $$\frac{1}{R^p}\mitQ{Q_{2R,2\rg}}{|u-u_{2R}|^p}
\le C\mitQ{Q_{2R,2\rg}}{|Du|^p} + CS^{(2-p)p}\left(\mitQ{Q_{2R,2\rg}}{|Du|^p})\right)^{p-1}.$$
By \mref{S1}, the right hand side is bounded by $CS^p$. Thus,
$$\frac{1}{R^p}\itQ{Q_{2R,2\rg}}{|u-u_{2R}|^p}= \frac{S^{2-p}R^{n+2}}{R^p}\mitQ{Q_{2R,2\rg}}{|u-u_{2R}|^p}
\le CR^{n+2}S^2.$$ 

\mref{S1} also gives 
$$\frac{1}{\rg}\itQ{Q_{2R,2\rg}}{|u-u_{2R}|^2} =
\frac{|Q_{2R,2\rg}|}{\rg}\mitQ{Q_{2R,2\rg}}{|u-u_{2R}|^2} \le
C|Q_{R,\rg}|S^p = CR^{n+2}S^2.$$

The above estimates and
\mref{uestpless2} yield $$\mypar{\sup_{t\in(-\rg,0)}\int_{B_{R}}|u-u_{x_0,R}^\chi(t)|^2dx}^{1-q/2}
\le C(R^{n+2}S^2)^{\frac{2}{n+2}} = CR^2S^{\frac{4}{n+2}}.$$

Hence, $$\barr{lll}\frac{1}{\rg}\mitQ{Q_{R,\rg}}{|u-u_{R}|^2} &\le&
CS^{p-2+\frac{4}{n+2}}\mitQ{Q_{R,\rg}}{|Du|^q} =
CS^{p-q}\mitQ{Q_{R,\rg}}{|Du|^q} \\&\le& \eg
S^p+C(\eg)\mypar{\mitQ{Q_{R,\rg}}{|Du|^q}}^{p/q}.\earr$$

From the assumption on $S$ and  the Caccioppoli inequality \mref{rev-cacc}
we derive the desired reverse
H\"older inequality in this case. \eproof

Finally, in order to obtain certain uniform continuity for the
integrals of gradients of weak solutions, we will need the following
measure theoretic result which could be of interest in its own right.

\blemm{revholder} Let $F,G_k$ ($k=1,\ldots,M$) be integrable
functions defined on $Q_{1,1}$ and $\ag,\bg,m_k$ be real numbers
with $\ag+1,\bg>0$ and $m_k\in(0,1)$. Assume that for any scaled
cylinders $Q_{R,\rg}\subset Q_{2R,2\rg}\subset Q_{1,1}$ with $\rg=S^{\ag}R^\bg$ and
$S\sim \mitQ{Q_{R,\rg}}{F}$ the following holds.

\beqno{w1} \mitQ{Q_{R,\rg}}{F} \le \eg \mitQ{Q_{2R,2\rg}}{F} +
\sum_{k=1}^{M}\mypar{\mitQ{Q_{2R,2\rg}}{G_k}}^{1/m_k}.\eeq

If $\eg>0$ is sufficiently small (depending on $
\|F\|_{L^1(Q_{1,1})}, \ag,\bg,m_k$)  then for any subset $A$ of
$Q_{\frac12,\frac12}$, $m\in(0,1)$ and  \beqno{Tkey} t^{\ag+1} \ge
\mitQ{Q_{1,1}}{F},\eeq 
there is a positive constant $C=C(n)$ such
that
$$\itQ{\Fg_t}{F} \le Ct^{1-m}\itQ{\Fg_t}{F^m} +
C\sum_{k=1}^Mt^{1-m_k}\itQ{\Gamma^k_t}{G_k}.$$ Here, for any $t>0$ ,
we set
 $\Fg_t = \{z\,:\, z\in A \mbox{ and } F(z)>t\}$ and $\Gamma^k_t =
\{z\,:\, z\in A \mbox{ and } G_k(z)>t^{m_k}\}$.

\elemm

\bproof For simplicity, we will consider the case when $M=1$ since
it is easy to extend the argument to the case $M>1$.

Let $P=Q_{\frac12,\frac12}$ and $t=\llg_0S$ with $\llg_0=\llg_0(n)$
being a constant to be determined later. We have from \mref{Tkey} \beqno{Skey}
\llg_0^{\ag+1}S^{\ag+1} \ge \mitQ{Q_{1,1}}{F}.\eeq

Define \beqno{Jdef}J(r) = \mitQ{Q_{r,S^\ag r^\bg}}{F}, \quad r>0.\eeq 

Obviously, for $r\in[\frac{1}{10},\frac12]$ and
$\llg_0=\llg_0(n)$ sufficiently small  we can have
\beqno{Schoice} J(r)\le S^{-\ag}C(n)\mitQ{Q_{1,1}}{F}=SC(n)\frac{1}{S^{\ag+1}}\mitQ{Q_{1,1}}{F}\le
\frac14S \quad \forall r\in[\frac{1}{10},\frac12].\eeq

Consider a point $z\in P$ such that $F(z) > S$.
Lebesgue's theorem yields $\lim_{r\to0}J(r)>S$. Thus, by the continuity of
the integral and the above inequality, we can find $r(z)\in(0,\frac{1}{10})$ such that
$J(r(z))=\frac12S$ and $J(\rg)\le\frac{1}{2}S$ for any
$\rg\in[r(z),\frac{1}{10}]$. This and \mref{Schoice} imply
that $J(2r(z))$ and $J(5r(z))$ are bounded by $\frac12 S$. Moreover,
there is a constant $c_0$ depending on $n$ such that $J(r(z))\le
c_0J(2r(z))$ and therefore $J(2r(z))\ge \frac12c_0^{-1}S$.

Hence, for $\rg(z)=S^{\ag}r^\bg(z)$  \beqno{cps} S \sim
\mitQ{Q_{2r(z),S^{\ag}(2r)^\bg(z)}}{F}\mbox{ and }
\mitQ{Q_{5r(z),S^{\ag}(5r)^\bg(z)}}{F} \le \frac12S.\eeq

We apply the Calderon-Zygmund lemma to $P$ to obtain a countable
family of disjoint subcubes
$\{Q_i\}=\{Q_{2r(z_i),S^{\ag}(2r)^\bg(z_i)}\}$ such that 
$\{z\,:\, z\in P\mbox{ and } F(z)>S\}\subset \cup_i \hat{Q}_i$. Here, $\hat{Q}_{i} =
Q_{5r(z_i),S^{\ag}(5r)^\bg(z_i)}$.

Let $A$ be a subset of $P$ and $\Fg_t = \{z\,:\, z\in A \mbox{ and }
F(z)>t\}$ and $\Gamma_t = \{z\,:\, z\in A \mbox{ and }
G(z)>t^{m_1}\}$. We will only consider subcubes $Q_i$'s such that
$Q_i\cap A\ne\emptyset$.

For such a subcube, $Q_{2r(z_i),S^{\ag}(2r)^\bg(z_i)}\subset Q_{4r(z_i),S^{\ag}(4r)^\bg(z_i)}\subset Q_{1,1}$ and by \mref{cps} we see that \mref{w1} holds and we
have two cases (with $\bar{Q}_{r,\rg} = Q_{2r,2\rg}$)
\beqno{w3}\mitQ{Q_i}{F} \le 2\eg \mitQ{\bar{Q}_i}{F} \mbox{ or }
\mitQ{Q_i}{F} \le 2\mypar{\mitQ{\bar{Q}_i}{G}}^{1/{m_1}}.\eeq

If the second case of \mref{w3} holds then because $S \sim
\mitQ{Q_i}{F}$ we have  the following
$$c_1(n)^{-m_1}S^{m_1} < \mypar{\mitQ{Q_i}{F}}^{m_1} \le 2^{m_1}\mitQ{\bar{Q}_i}{G}\le
c_2(n)2^{m_1}\mitQ{\hat{Q}_i}{G}.$$ Hence, by splitting the integral on $\hat{Q}_i$ into those on $\hat{Q}_i\cap \Gamma_t$ and $\hat{Q}_i\setminus \Gamma_t$, we have for some $c=c(n)$
$$S^{m_1} |\hat{Q}_i| < (2c)^{m_1}\itQ{\hat{Q}_i\cap \Gamma_t}{G} + (2c)^{m_1}
\itQ{\hat{Q}_i\setminus \Gamma_t}{G}.$$ This gives
$$S^{m_1}|\hat{Q}_i| <
(2c)^{m_1}\itQ{\hat{Q}_i\cap \Gamma_t}{G} + (2c)^{m_1}
t^{m_1}|\hat{Q}_i|$$ and furthermore $$S^{m_1}|\hat{Q}_i| \le
c_3(n)\itQ{\hat{Q}_i\cap \Gamma_t}{G}, $$ if $\llg_0=t/S$ is
sufficiently small (such that $(2c)^{m_1}t^{m_1}\le\frac12S^{m_1}$ and \mref{Schoice} still holds). We then
fix such $\llg_0$.

Arguing similarly, with $G,m_1$ now being $\eg F$ and $1$, we see that if
the first case of \mref{w3} holds and $\eg$ is small then
$$S|\hat{Q}_i| \le 4\eg\itQ{\hat{Q}_i\cap \Fg_S}{F}.$$

Together, we have $$ S|\hat{Q}_i| \le 4\eg \itQ{\hat{Q}_i\cap
\Fg_S}{F} + c_3(n)S^{1-m_1}\itQ{\hat{Q}_i\cap \Gamma_t}{G}.$$

Since $\Fg_S\subset\cup_i\hat{Q}_i$, we
now have by \mref{cps} and the fact that $Q_i$'s are disjoint
$$\itQ{\Fg_S}{F} \le\sum_{\hat{Q}_i\cap\Fg_S\ne\emptyset}
\itQ{\hat{Q}_i}{F} \le
\frac12S\sum_{\hat{Q}_i\cap\Fg_S\ne\emptyset}|\hat{Q}_i|
=C(n)S\sum_i|Q_i| = C(n)S\left|\cup {Q}_i\right|.$$ By Vitali's
covering lemma, we can find a subsequence of disjoint subcubes
$\{\Pi_i\}$ of $\{\hat{Q}_i\}$ such that $\cup{Q}_i\subset
\cup\hat{\Pi}_i$ and therefore
$$\left|\cup {Q}_i\right| \le \left|\cup \hat{\Pi}_i\right| \le
\sum |\hat{\Pi}_i| \le C_1(n) \sum |\Pi_i|.$$

Thus, as $\Pi$'s are disjoint, we have from the above estimates for $\Pi_i=\hat{Q}_i$ that
$$\barr{lll}\itQ{\Fg_S}{F} &\le& C_2(n)S \sum |\Pi_i| \le
C_3(n)\sum\mypar{4\eg \itQ{\Pi_i\cap \Fg_S}{F} +
S^{1-m_1}\itQ{\Pi_i\cap \Gamma_t}{G}}\\ &\le& C_3(n)4\eg \itQ{
\Fg_S}{F} + C_3(n)S^{1-m_1}\itQ{\Gamma_t}{G}.\earr$$

On the other hand, as $\llg_0$ is small, $t<S$ and therefore $\Fg_S\subset \Fg_t$. We then have
$$\itQ{\Fg_t\setminus \Fg_S}{F}\le S^{1-m} \itQ{\Fg_t}{F^m} \le 2C_4(n)
t^{1-m}\itQ{\Fg_t}{F^m}.$$

Hence, by choosing $\eg$ sufficiently small, we get $$\itQ{\Fg_t}{F}
\le Ct^{1-m}\itQ{\Fg_t}{F^m} + Ct^{1-m_1}\itQ{\Gamma_t}{G}.$$

This completes the proof. \eproof

\brem{Srem} Note that the number $S$ in the proof is fixed and needs
only satisfy \mref{Skey}. \erem

We then have the following result on the uniform continuity of
integrals.

\blemm{weakconv} Let $\{F^{(i)}\}$ and $\{G_k^{(i)}\}$,
$k=1,\ldots,M$, be bounded sequences of functions in $L^1(Q_{1,1})$
satisfying \mref{w1} of \reflemm{revholder}. Assume that the
sequences $\{G_k^{(i)}\}$ are weakly convergent in $L^1(Q_{1,1})$. If
$\eg$ is sufficiently small then the integrals of $F^{(i)}$ are uniformly continuous in the following sense: for any given $\dg>0$ there is
$\mu(\dg)>0$ such that $$ \itQ{A}{F^{(i)}} < \dg \mbox{ for all $i$
if $A\subset Q_{\frac12,\frac12}$ and $|A|\le \mu(\dg)$.}$$ \elemm

\bproof We can again consider the case $M=1$. Fix a $t$ satisfying
\mref{Tkey} of \reflemm{revholder}. We have shown that if $\eg$ is
small enough then we have for $F=F^{(i)}$ and $G=G^{(i)}$.
$$\barr{lll}\itQ{A}{F} &=& \itQ{A\setminus\Fg_t}{F} +
\itQ{\Fg_t}{F} \le t|A|+ Ct^{1-m}\itQ{\Fg_t}{F^m} +
Ct^{1-m_1}\itQ{\Gamma_t}{G}\\&\le& t|A|+ Ct^{1-m}\itQ{A}{F^m} +
Ct^{1-m_1}\itQ{A}{G}.\earr$$

Since $|F^{(i)}|^m, G^{(i)}$ are weakly convergent in $L^1(Q_{1,1})$, we
can apply \cite[Corollary IV.11]{DF} on the uniform continuity of
integrals to see that the last two integrals are uniformly small if
$|A|$ is small. The assertion then follows easily. \eproof

Another consequence of \reflemm{revholder} is the following higher
integrability result.

\blemm{higherLq} Let $\{F^{(i)}\}$ and $\{G_k^{(i)}\}$,
$k=1,\ldots,M$, be bounded sequences of functions in $L^1(Q_{1,1})$
satisfying \mref{w1} of \reflemm{revholder} (with $\rg=S^{\ag}R^\bg$).  If $\eg$ is
sufficiently small then  there is  some $r>1$ such that if $G_i\in
L^{r-m_i+1}(Q_{1,1})$ then the following estimate  holds \beqno{Lqest}
\itQ{Q_{\frac12,\frac12}}{F^r} \le
C\left(\itQ{Q_{1,1}}{F}\right)^{\frac{r-1}{1+\ag}+1}+
C\Sg\itQ{Q_{1,1}}{G_i^{r-m_i+1}}.\eeq \elemm

\bproof For $m\in(0,1)$ and $r>1$, we define $$\fg(t)=\itQ{\Fg_t}{F^m},\quad
\og_i(t)=\itQ{\Gamma_t}{G_i},$$ and $$I_r(t) = \itQ{\Fg_t}{F^r}.$$

For $A=Q_{\frac12,\frac12}$, the assertion in \reflemm{revholder} can be
written as
$$-\int_t^\infty\tau^{1-m}d\fg(\tau) \le C[t^{1-m}\fg(t) + \Sg_i
t^{1-m_i}\og_i(t)] \quad \forall t\ge a:=\left(\mitQ{Q_{1,1}}{F}\right)^\frac1{\ag+1}.$$

A simple modification of the Gehring lemma in \cite[Lemma 6.3,
p.200]{Giusti} provides some $r>1$ such that $$-\int_a^\infty
u^{r-m}d\fg(u) \le -Ca^{r-1}\int_a^\infty u^{1-m}d\fg(u) -
C\Sg\int_a^\infty u^{r-m_i}d\og_i(u).$$

This gives (see \cite{Giusti})
$$\itQ{Q_{\frac12,\frac12}}{F^r}
\le Ca^{r-1}\itQ{Q_{1,1}}{F}+ C\Sg\itQ{Q_{1,1}}{G_i^{r-m_i+1}}.$$ The
definition of $a$ then gives the desired \mref{Lqest}. \eproof

We now consider a sequence of vector functions $u_k$ which almost
solve \mref{plemeqn} in the following sense: \beqno{plemeqn1} \left|
\itQ{Q_{R,\rg}}{-u\fg_t+\myprod{A_k(u,Du),D\fg}}\right|\le\dg
\|Du\|^{\frac pq}_{L^p(Q_{R,\rg})}\|D\fg\|_{L^p(Q_{R,\rg})}\eeq for all $\fg\in
V_p^0(Q_{R,\rg})$, $Q_{R,\rg}\subset Q_{1,1}$.

\blemm{duweakconv} Assume that $2n/(n+2)<p<2$. Let $\{u_k\}$ be a
sequence of vector functions satisfying \mref{plemeqn1} and the
norms $\|Du_k\|_{L^p(Q_{1,1})}$ are uniformly bounded. The sequence of
matrices $A_k$ is assumed to satisfied the ellipticity condition E)
of \refsec{mainsec}. If $\dg$ is sufficiently small then the integrals of  $|Du_k|^p$ are uniformly continuous. \elemm

\bproof  We will apply \reflemm{weakconv} here by taking
$F^{(i)}=|Du_i|^p$ and $G_k^{(i)}$ ($k=1,2$) to be either $|Du_i|^q$
or $|Du_i|^{p-1}$ with $q=np/(n+p)$. Consider first the case $p>2$.
It is easy to see that $u_i$ satisfies the Poincar\'e
and Caccioppoli type inequalities of \reflemm{pineq} and
\reflemm{CACCineq}. Hence the reverse H\"older inequality
\mref{revlp} of \reflemm{reverse2} holds for $Du_i$ with $\eg=2\dg$
and $$  \rg=S^{2-p}R^2 \mbox{ and } S^p\sim \mitQ{Q_{R,\rg}}{|Du|^p}.$$  Hence,
 the assumption
\mref{w1} of \reflemm{revholder} is verified with $\ag,S$ there being $(2-p)/p,S^p$ respectively, $\rg=S^{\ag}R^\bg$ and $\bg=2$. 

As the norms $\|Du_k\|_{L^p(Q_{1,1})}$ are uniformly bounded, $G_i$'s are
uniformly bounded in $L^r$ for some $r>1$ and they are weakly
convergent in $L^1$. \reflemm{weakconv} then applies here to give
our lemma if $\eg$, or $\dg$, is sufficiently small.

If $p<2$, as
in the proof of \reflemm{revholder}, we fix a number $S$ such that (see \mref{Skey})$$\llg_0^\frac{2}{p}S^2 \ge \left(
\mitQ{Q_{1,1}}{|Du_k|^p} + \frac{1}{S^{2-p}}\mitQ{Q_{1,1}}{|u_k-(u_k)_{Q_{1,1}}|^2}\right)$$ 

Define (see \mref{Jdef}) $$J(r) = 
\mitQ{Q_{r,S^{2-p}r^2}}{|Du_k|^p} + \frac{1}{S^{2-p}r^2}\mitQ{Q_{r,S^{2-p}r^2}}{|u_k-(u_k)_{Q_{1,1}}|^2}.$$ 

We will be interested in the
set where $|Du_k|^p>S^p$. At each point $z$ of this set, the argument leading to \mref{cps} in the proof allows us to
find a cylinder $Q_{R,S^{2-p}R^2}(z)$ and a positive constant $c_1(n)$ such that $$c_1(n)\left( \mitQ{Q_{R,S^{2-p}R^2}}{|Du|^p} +
\frac{1}{S^{2-p}R^2}\mitQ{Q_{R,S^{2-p}R^2}}{|u-u_R|^2}\right)\le S^p.$$ Therefore, the condition
\mref{S1} on $S$ is verified and a reverse H\"older inequality
\mref{w1} for $F^{(k)}=|Du_k|^p$ is available again. Noting that $S$ is fixed, we see that the proof can continue
as before. \eproof

We also have the following $L^q$ estimates for $Du$ as a result of
\reflemm{higherLq}.

\blemm{duLq} Assume that $p>2n/(n+2)$. Consider a vector functions
$u$ satisfying \mref{plemeqn1} with $A_k$ is assumed to satisfied
the ellipticity condition E) of \refsec{mainsec}. If $\dg$ is
sufficiently small then there is $\eg>0$ such that
\beqno{lqdu}R^{p+\eg}\mitQ{Q_{\frac12R,\frac12R^p}}{|Du|^{p+\eg}}
\le
C\eeq for some constant $C$ depending on $R^p\mitQ{Q_{R,R^p}}{|Du|^p}$. If $2n/(n+2)<p<2$, the above constant $C$ also depends on $\frac{1}{R^p}\mitQ{Q_{R,R^p}}{|u-u_R|^2}$.

\elemm

\bproof We make the scaling $x\to x/R$, $t\to t/R^p$ in the equation for
$u$ and need only show that \mref{lqdu} when $R=1$. We now set $F=|Du|^p$ and $G_i$ be
either $|Du|^{p-1}$ or $|Du|^q$ with $q=np/(n+p)$. Accordingly,
$m_i$ will be either $(p-1)/p$ or $q/p=n/(n+p)$ respectively. Thus,
$G_i=F^{m_i}$ and belongs to $L^{r-m_i+1}$ if $(r-m_i+1)m_i\le1$. We
can find such $r>1$ if  $(2-m_i)m_i<1$ but this requirement is just one of the followings
$$(2-\frac{p-1}{p})\frac{p-1}{p} < 1 \Leftrightarrow p^2-1<p^2
\mbox{ and } (2-\frac{n}{n+p})\frac{n}{n+p} < 1 \Leftrightarrow
n(n+2p)<(n+p)^2.$$

Thus, \reflemm{higherLq} applies here with $\ag=(2-p)/p$ and
$r=1+\eg/p$ to give \mref{lqdu} when $R=1$. The dependence of the constant $C$ on $\frac{1}{R^p}\mitQ{Q_{R,R^p}}{|u-u_R|^2}$ comes from the choice of $S$ in the proof of \reflemm{duweakconv}. \eproof

\brem{duLq-rem} When $p>2$, it is easy to see that the quantity $$\left(R^{p+\eg}\mitQ{Q_{\frac12R,\frac12R^p}}{|Du|^{p+\eg}}\right)^\frac{p}{p+\eg} $$ can be bounded by $$
C_0\max\{\left(R^p\mitQ{Q_{R,R^p}}{|Du|^{p}}\right)^\frac{p(2+\eg)}{2(p+\eg)},R^p\mitQ{Q_{R,R^p}}{|Du|^{p}}\}$$ for some constant $C_0$ independent of $Du$. We also remark that the exponent $\frac{p(2+\eg)}{2(p+\eg)}>1$ when $p>2$.
\erem

\section{The approximation lemmas}\eqnoset\label{app-sec}

\newc{\q}{Q_1}
\newc{\qq}{Q_{\frac12}}
\newc{\qqq}{Q_{\frac23}}
\newc{\sq}{S_1}
\newc{\sqq}{S_{\frac12}}
\newc{\sqqq}{S_{\frac23}}

In this section, we establish one of the main tools of our work -
the nonlinear approximation lemma for $p$-Laplacian systems. To
begin, let us fix a cylinder $Q_R=B_R\times[-R^p,0]$ and consider
two collections $\ccA$ of matrix-valued functions $A$ and $\ccB$ of
vector valued functions $B$ satisfying the followings

\bdes \item[a.1)] There are positive constants $\llg,\LLg$ such that
the ellipticity condition \mref{aaa}
in E) holds for each $A\in\ccA$ .

\item[a.2)] For any $B\in\ccB$, $B(u,\zeta)$ is Lipschitz in $u$ and there
is a positive constant $C$ such that \beqno{pfff} |B(u,\zeta)| \le C +
C|\zeta|^{p-1}.\eeq  
\item[a.3)] For each $A\in\ccA$, $B\in\ccB$ and any given function $g\in
C^1(Q_{\rg})$, where $Q_\rg\subset Q_R$, the system
$$\itQ{Q_\rg}{-u\fg_t+\myprod{A(u,Du),D\fg}-\myprod{B(u,Du),\fg}}=0, \,\forall \fg\in C^1_0(Q_{\rg}), \quad u=g \mbox{ on $S_{\rg}$}
$$ has a bounded weak solution $u$ with
$\|u\|_{L^\infty(Q_\rg)} \le C(\|g\|_{L^\infty(Q_\rg)})$.

\item[a.4)] (Monotonicity) There is a positive constant $\llg_0$  such that
$$ \itQ{Q_\rg}{\langle A(u,Du)-A(v,Dv),Du-Dv\rangle}\ge
\llg_0\itQ{Q_\rg}{|Du-Dv)|^p}$$ for any $u,v\in V_p(Q_\rg)$.

\edes

The monotonicity condition a.4) has been frequently assumed in
literature concerning the uniqueness of weak solutions to the systems
described in a.3). We will state the first version of our nonlinear
approximation results under this condition in order to streamline our presentation and ideas.
Later, we will replace a.4) by more practical assumptions and the
proof will be similar modulo some technical modifications.

We will first prove that

\bprop{heat-aa} Assume a.1)- a.4) and $p>2n/(n+2)$. For any given $M,\eg>0$
and $\bg>1$ there exists $\dg\in(0,1]$ that depends only on
$\llg,\LLg,M,\eg$ such that if $A\in\ccA, B\in\ccB$ and
 $u\in V_p(Q_R)$ satisfying

\beqno{ha1a} \mitQ{Q_R}{|u|^2} +R^p\mitQ{Q_R}{|Du|^p}\le M,\eeq
and\beqno{ha1aa} \left| \itQ{Q_R}{-u\fg_t+\myprod{A(u,Du),D\fg}
-\myprod{B(u,Du),\fg} }\right|\le\dg \|Du\|^{\frac
pq}_{L^p(Q_R)}\|D\fg\|_{L^p(Q_R)}\eeq for all $\fg\in V_p^0(Q_R)$,
then either \beqno{hb3a} R^p\mitQ{Q_R}{|Du|^p} < \eg, \eeq or there
exists $v\in V_p(Q_{R/2})$ such that
$$\mitQ{Q_{R/2}}{|v|^2} +R^p\mitQ{Q_{R/2}}{|Dv|^p}\le C(\llg,\LLg,\llg_0,M)$$ and
\beqno{ha3aa}\itQ{Q_{R/2}}{-v\fg_t+\myprod{A(v,Dv),D\fg}
-\myprod{B(v,Dv),\fg} }=0\eeq for all $\fg\in C_0^1(Q_{R/2})$, and
\beqno{ha3a}
\left\{\barr{l}\mitQ{Q_{R/2}}{|v-u|^2}+\mitQ{Q_{R/2}}{|v-u|^p} \le
\eg R^{p}\mitQ{Q_R}{|Du|^p},\\  \mitQ{Q_{R/2}}{|Dv|^p}\le
\bg\mitQ{Q_{R/2}}{|Du|^p}.\earr\right.\eeq
 \eprop

\bproof For simplicity we will present the proof when $B$ is
identically zero. It is not difficult to see that the presence of
$B$, satisfying our assumptions, will introduce some extra terms
which can be easily treated by the same argument and a simple use of Young's inequality.

The proof is by contradiction, we then assume that  there exist
$\eg_0>0$ and sequences $\{A_k\},\{u_k\}, \{\eg_k\}$ and cylinders
$Q_{R_k}(x_k,t_k)$ such that for $Q_{R_k}=Q_{R_k,R_k^p}(x_k,t_k)$ we have
\beqno{ukeqn} |\itQ{Q_{R_k}}{-u_k\fg_t+\myprod{A_k(u_k,Du_k),D\fg}}| \le
\eg_k\|Du_k\|_{L^p(Q_{R_k})}^\frac pq\|D\fg\|_{L^q(Q_{R_k})}\eeq for all
$\fg\in V_p^0(Q_{R_k})$ but \beqno{Du-alt1}
\limsup_{k\to\infty}R_k^p\mitQ{Q_{R_k}}{|Du_k|^p} >0,\eeq and \beqno{ukeqn1}
\mitQ{Q_{R_k/2}}{|v-u_k|^2}+\mitQ{Q_{R_k/2}}{|v-u_k|^p} >
\eg_0 R_k^{p}\mitQ{Q_{R_k}}{|Du_k|^p}\eeq for all
$v$ satisfying $$ \itQ{Q_{R_k}}{-v\fg_t+\myprod{A_k(v,Dv),D\fg}}=0 \quad
\mbox{ for all $\fg\in V_p^0(Q_{R_k/2})$}.$$

We then make a change of variables
$$\widetilde{u}_k(x,t) = u_k(x_k+ R_k x, t_k+ R^p_k t), \quad (x,t)\in
Q_1.$$ By the boundedness assumption \mref{ha1a}, the norms
$\|u_k\|_{V_p(Q_{R_k})}$ are uniformly bounded by $M$ and so are
$\|\widetilde{u}_k\|_{V_p(Q_1)}$. Thus, by scaling and translation, we can assume
$R=1$ in \mref{ukeqn}-\mref{ukeqn1} and note that \mref{ha1a} and
\mref{ha1aa} remain. This also proves that $\dg$ is independent of
$R$.

For any positive real number $h$ and any vector valued function $f$ in $L^1(Q_1)$,
we denote by $f_h=J_h*f$ the standard mollifier of $f$. That is, for
some smooth nonnegative function $J$ with compact support in the
unit ball $Q_1$ of $\RR^{n+1}$ and $\|J\|_{L^1(\RR^{n+1})}=1$, we
write
$$f_h(Z)=J_h*f(Z)=\frac{1}{h^{n+1}}\itQ{Q_1}{J\left(\frac{|Z-z|}{h}\right)f(z)}, \quad Z\in\RR^{n+1}.$$

 Let $\{h_{k}\}$ be some sequence
of positive reals converges to $0$ and $g_{k}=(u_k)_{h_{k}}$, which
is in $C^1(\q)$. We then define $U_k$ to be the solutions of

\beqno{asys12b}\left\{\barr{ll}
\itQ{\qqq}{-U_k\fg_t+\myprod{A_k(U_k,DU_k),D\fg}}=0 & \forall \fg\in
C^1_0(\qqq)\\  U_k=g_{{k}}& \mbox{ on $\sqqq$.}\earr\right.\eeq

Note that, by a.3), $\|U_k\|_\infty\le C(\|g_k\|_\infty)$.

The following claims provide a contradiction to \mref{Du-alt1} and
\mref{ukeqn1} and prove our proposition.

\bdes\item[Claim I:] There is a constant $C$ such that
$\|DU_k\|_{L^p(\qqq)}\le C\|Du_k\|_{L^p(\qqq)}$.
\item[Claim II:]  $u_k-U_k\to 0$ in $L^2(\qqq)$ and $L^p(\qqq)$.\edes

{\bf Proof of Claim I:} Let $(\ag)_k=a_k(u_k,Du_k)_{h_k}$ and replace
$\fg$ in the inequality for $u_k$ by $\fg_{-h_k}$, whose support is
in $Q_{5/6}$ if $h_k$ is sufficiently small. From \mref{ukeqn}, the
following holds

$$ |\itQ{\q}{[-\myprod{g_k,\fg_t} + \myprod{(\ag)_k, D\fg}] }|\le
\eg_k\|Du_k\|^{\frac pq}_{L^p(\q)}\|D\fg\|_{L^p(\q)},
 \forall \fg\in C^1_0(\q).$$

For any  $\tau\in(-(\frac23)^p,0)$ and sufficiently small positive $h$, let $\eta(s)$ be a $C^1$ function such that $\eta(s)=0$ if $s<\tau$, $\eta(s)=0$ if $s>\tau+h$ and $\eta$ is almost linear in $(\tau,\tau+h)$. Subtracting this inequality by the equation of $U_k$ and testing the
result with $(U_k(x,s)-g_k(x,s))\eta^2(s)$, which vanishes on the boundary of $\qqq$, we easily
obtain the following for $\fg=U_k-g_k$ when we send $h$ to 0
$$\sup_{s\in[-\frac49,\tau]}\iidx{B_\frac23}{|\fg|^2}+\itQ{\qqq^\tau}{\myprod{A_k(U_k,DU_k)-(\ag)_k,D\fg}} \le
\eg_k\|Du_k\|^{\frac pq}_{L^p(\q)}\|D\fg\|_{L^p(\q)}.$$ 
Here, for
any cylinder $Q$, we denote $Q^t=Q\cap \{(x,s)\,:\, s\le t\}$. The above then implies
$$ \barr{ll}\lefteqn{\sup_{s\in[-\frac49,\tau]}\iidx{B_\frac23}{|\fg|^2} +
\itQ{\qqq^\tau}{\myprod{a(U_k,DU_k),DU_k}}=}\hspace{.5cm}&\\
& \itQ{\qqq^\tau}{\myprod{A_k(U_k,DU_k),Dg_k} + \myprod{(\ag)_{k},
(DU_k-Dg_k)}} +\eg_k\|Du_k\|^{\frac
pq}_{L^p(\q)}\|D\fg_{h_k}\|_{L^p(\q)}.\earr$$

Since $\|D(U_k)_{h_k}\|_{L^p(\q)}\le C\|DU_k\|_{L^p(\q)}$,  a simple
use of the Young inequality with $\eg_k$ sufficiently small and the
ellipticity of $A_k$ in a.1) give \beqno{ukgk}
\barr{ll}\lefteqn{\sup_{s\in[-\frac49,\tau]}\iidx{B_\frac23}{|U_k-g_k|^2}
+ \itQ{\qqq^\tau}{|DU_k|^p}\le}\hspace{1cm}&\\ &C\itQ{\qqq^\tau}{[|Dg_k|^p
+ |(\ag)_{k}|^q +|Du_k|^p]}+\eg_k\|Du_k\|^{\frac
pq}_{L^p(\q)}\|D\fg_{h_k}\|_{L^p(\q)}.\earr\eeq

Because $\|Dg_k\|_{L^p(\qqq^\tau)}\le C\|Du_k\|_{L^p(\qqq^\tau)}$ and
$\|(\ag)_{k}\|^q_{L^q(\qqq^\tau)}=
\|A_k(u_k,Du_k)_{h_k}\|^q_{L^q(\qqq^\tau)} \le
C\|Du_k\|^p_{L^p(\qqq^\tau)}$, we obtain from the above and the Young inequality the following estimate
\beqno{dukest}\|DU_k\|^p_{L^p(\qqq^\tau)} \le C\|Du_k\|^p_{L^p(\qqq^\tau)}
+ \eg_k\|DU_k\|^p_{L^p(\qqq)}\quad \forall \tau\in[-(\frac23)^p,0].\eeq
This established our first claim if we take $\tau=0$.

{\bf Proof of Claim II:} Now, for any $\rg,r\in(0,\frac23)$, we
write $Q'_1=\qqq\cap\{(x,s):\, s\le -(\frac23)^p+\rg^2$,
$Q'_2=\qqq\cap \{(x,t):\, |x|>\frac23-r\}$ and $Q'_3=Q'_1\cap Q'_2$.
These sets are the thin layers at the base and lateral sides of the
cylinder $\qqq$.

Let  $\fg(x,t)=\psi(x)\eta(t)$, where $\psi,\eta$
are respectively cut-off functions in $x,t$. That is, $\psi$ is a
cut-off function for $B_{\frac23-r}$ and $B_\frac23$ and $\eta$ is a
cut-off function for $[-(\frac23)^p+\rg^2,0]$ and
$[-(\frac23)^p,0]$. We can assume that $|D\psi|\le1/r$ and
$|\eta_t|\le 1/\rg^2$.

Denote $H_k=u_k-U_k$ and $\Fg=H_k\fg^2$. Testing the equations for $U_k$ with $\Fg$ and
replacing $\fg$ in \mref{ukeqn} by $\Fg$, we get by subtracting the two
results
$$ \left|\itQ{\qqq}{\left[H_k\pder{\Fg}{t}+
\myprod{A_k(u_k,Du_k)-A_k(U_k,DU_k),D\Fg}\right]}\right| \le
\eg_k\|Du_k\|^{p/q}_{L^p(\qqq)}\|D\Fg\|_{L^p(\qqq)}.$$

We now write $\bbA = A_k(u_k,Du_k)-A_k(U_k,DU_k)$ and derive from
the above the following
\beqno{bbA}\barr{lll}\lefteqn{\sup_{s\in[-(\frac23)^p,0]}
\iidx{B_\rg}{H_k^2\fg^2} + \itQ{\qqq}{\myprod{\bbA, D\Fg}}
\le}\hspace{2cm}&&\\ &&\itQ{Q}{|H_k|^2\fg\left|\pder{\fg}{t}\right|}
+ \eg_k\|Du_k\|^{p/q}_{L^p(\qqq)}\|D\Fg\|_{L^p(\qqq)}.\earr\eeq

We consider the second term on the left and note that $D\Fg=\fg^2
DH_k + 2H_k \fg D\fg$. By the monotonicity assumption a.4) with $u=u_k$ and $v=U_k$, we have
$$\bbA DH_k\fg^2 \ge \llg_0\fg^2|DH_k|^p.$$

On the other hand, because $D\fg=0$  in $Q\setminus Q'_2$,
$|D\fg|\le 1/r$  and $|\bbA| \le C[|Du_k|^{p-1}+|DU_k|^{p-1}]$, we
also have via the Young inequality
$$\barr{lll} \itQ{\qqq}{\bbA H_k\fg D\fg} &\le &  C\itQ{\qqq}{(
|Du_k|^{p-1}+|DU_k|^{p-1})|H_k||\fg| |D\fg|} \\
&\le& C\itQ{Q'_2}{[|Du_k|^p +  |DU_k|^p]} +
C\frac{1}{r^p}\itQ{Q'_2}{H_k^p}.\earr$$

Also, because $\fg_t=0$ in $Q\setminus Q'_1$ and $|\fg_t|\le
1/\rg^2$, we have
$$ \itQ{Q}{|H_k|^2\fg\left|\pder{\fg}{t}\right|} \le
C\frac{1}{\rg^2}\itQ{Q'_1}{|H_k|^2}. $$

On the other hand, as $$\|D\Fg\|^p_{L^p(\qqq)} \le
C\|DH_k\|^p_{L^p(\qqq)} + C\frac{1}{r^p}\itQ{Q'_2}{H_k^p},$$ a simple
use of the Young inequality and the above estimates allow us to
deduce from \mref{bbA} that

\beqno{key1}\sup_{s\in[-(\frac23)^p,0]}\iidx{B_\rg}{H_k^2\fg^2}
+\itQ{\qqq}{\llg_0 |DH_k|^p\fg^p}\le F_k,\eeq where
\beqno{key2}\barr{lll}F_k &=& \eg_k
\itQ{\qqq}{(|Du_k|^p+|DU_k|^p)}+C\itQ{Q'_2}{[|Du_k|^p +
|DU_k|^p]}+\\&&C\frac{1}{r^p}\itQ{Q'_2}{H_k^p}+C\frac{1}{\rg^2}\itQ{Q'_1}{H_k^2}.\earr\eeq

For any given $\eg>0$ we will show that if $r,\rg$ are sufficiently
small and $k$ is large then
\beqno{s1}\frac{1}{r^p}\itQ{Q'_2}{|H_k|^p}+\frac{1}{\rg^2}\itQ{Q'_1}{|H_k|^2}<\eg,
\eeq and \beqno{s2}\itQ{Q'_2}{(|Du_k|^p + |DU_k|^p)} <\eg.\eeq The
above estimates yield
$$ \sup_{s\in[-(\frac23)^p,0]} \iidx{B_\rg}{H_k^2} +\itQ{\qqq}{\llg_0 |DH_k|^p\fg^p}\le C\eg_k+2\eg,$$ and thus
$\|H_k\|_{L^2(\qq)}<\eg$. By Sobolev's imbedding inequality (see
\cite[Proposition 3.1, p.7]{Dib}), we also have
$\|H_k\|_{L^p(\qq)}<\eg$. Choosing $\eg$ small we obtain the desired
contradiction.

Concerning \mref{s1}, when $r,\rg$ have been fixed, because
$H_k=(u_k-g_k)+(g_k-U_k)$ and $u_k-g_k$ converge to $0$ in $L^p$ and
$L^2$, we need only to prove that \beqno{ugk}
\frac{1}{r^p}\itQ{Q'_2}{|g_k-U_k|^p}+\frac{1}{\rg^2}\itQ{Q'_1}{|g_k-U_k|^2}<\eg,\eeq
if $k$ is large and $r,\rg$ are sufficiently small (uniformly in
$k$).

For the integral on $Q'_1$, from \mref{ukgk}, we find that
$$\barr{ll}\lefteqn{\itQ{Q'_1}{|U_k-g_k|^2} \le C\rg^2\sup_{s\in[-(\frac23)^p,-(\frac23)^p+\rg^2]}
\iidx{B_\frac23,s}{|U_k-g_k|^2}\le}\hspace{2cm}&\\&
C\rg^2\itQ{Q'_1}{[|Dg_k|^p +
|(\ag)_{h_k}|^q]}+\eg_k\rg^2\|DU_k\|^p_{L^p(\qqq)}.\earr$$

On the other hand, because $U_k=g_k$ on the lateral part of $\sqqq$,
we can use the Poincar\'e inequality in $x$ to get
$$\frac{1}{r^p}\itQ{Q'_2}{|U_k-g_k|^p} \le \itQ{Q'_2}{|DU_k-Dg_k|^p}
\le \itQ{Q'_2}{|Dg_k|^p+|DU_k|^p}.$$

By \mref{ukgk} again, $$\barr{lll}\itQ{Q'_2}{|DU_k|^p}&\le&
\itQ{Q'_2\setminus Q'_1}{|DU_k|^p}+\itQ{Q'_1}{|DU_k|^p}\\&\le&
\itQ{Q'_2\setminus Q'_1}{|DU_k|^p}+\itQ{Q'_1}{[|Dg_k|^p +
|(\ag)_{h_k}|^q]}.\earr$$

But $$\itQ{Q'_1}{[|Dg_k|^p + |(\ag)_{h_k}|^q]} \le
C\itQ{Q'_1}{|Du_k|^p}.$$

Therefore, the left hand side of \mref{ugk} can be estimated by
\beqno{intsmall} C\left[\itQ{Q'_1}{|Du_k|^p} +\itQ{Q'_2}{|Du_k|^p}+\itQ{Q'_2\setminus
Q'_1}{|DU_k|^p}+\eg_k\right].\eeq

Obviously, the left hand side of \mref{s2} is also bounded by the
above. Thus, we need only prove that the integrals in
\mref{intsmall} can be arbitrarily small (uniformly in $k$) if
$r,\rg$ are sufficiently small. By the uniform continuity of integrals (see \reflemm{duweakconv}),  the integral of $|Du_k|^p$ over $Q'_1$ is small if
the measure $|Q'_1|$ , or $\rg$, is sufficiently small (but independent of $k$).
Hence, the first term of \mref{intsmall} is small.  Fixing
such a $\rg$, We then repeat the argument to see that if $r$ is
small then so is the second integral in \mref{intsmall}. Similarly, the
integral of $|DU_k|^p$ over $Q'_2\setminus Q'_1$ is small.

Therefore, the right hand side $F_k$ of \mref{key1} can be
arbitrarily small if $r,\rg$ are sufficiently small and $k$ is
large. As we mentioned earlier this gives the proof of the second
claim and completes our proof. \eproof

We now consider the following alternative of the monotonicity
condition a.4).

\bdes \item[a.4')] For any $w\in \RR^m$ and $U,V\in\RR^{nm}$, there
holds \beqno{AAmono} \myprod{A(w,U)-A(w,V),U-V} \ge
\llg_0|U-V|^2\left\{\barr{ll}\min\{|U|^{p-2},|V|^{p-2}\}&
U\ne0\mbox{ or } V\ne0\\0&\mbox{otherwise.}\earr\right. \eeq
Moreover, $A(u,\zeta)$ is Lipschitz in $u$ in the following sense
$$|A(u,\zeta)-A(v,\zeta)| \le C|u-v||\zeta|^{p-1}.$$\edes

Concerning the condition \mref{AAmono}, if $A(u,U)$ is differentiable in $U$ then we note that
$$\myprod{A(w,U)-A(w,V),U-V} =
\int_0^1\myprod{\pder{A}{\zeta}(w,sU+(1-s)V)ds(U-V),U-V}.$$
Therefore, \mref{AAmono} can be verified if the matrix
$\pder{A}{\zeta}$ is positive definite and $$
\myprod{\pder{A}{\zeta}(w,sU+(1-s)V)\eta,\eta}\ge
\llg(sU+(1-s)V)|\eta|^2 \quad \forall s\in(0,1),\eta\in\RR^{nm}$$
for some $\llg(sU+(1-s)V)\ge |sU+(1-s)V|^{p-2}$.

\bprop{heat-aac} The conclusion of \refprop{heat-aa} holds if the
monotonicity condition a.4) is replaced by a.4').\eprop

\bproof We revisit the proof of \refprop{heat-aa} and point out
necessary modifications under a.4'). As before, for
 $H_k=u_k-U_k$ and $\fg(x,t)=\psi(x)\eta(t)$
with $\psi,\eta$ being respectively cut-off functions in $x,t$ for
$\qqq$. That is, $\psi$ is a cut-off function for $B_{\frac23-r}$
and $B_\frac23$ and $\eta$ is a cut-off function for
$[-(\frac23)^p+\rg^2,\tau]$ and $[-(\frac23)^p,\tau]$ , where $\tau$
is any number in $[-(\frac23)^p,0]$.

Again, the proof is by contradiction and we can see that the proof
of claim I in the proof of \refprop{heat-aa} is still applicable here. We need only consider claim II.
First of all, the assumptions \mref{Du-alt1} and \mref{ukeqn1} give
that
$$\itQ{\qqq}{|H_k|^2} \ge \eg_1>0$$ for some fixed $\eg_1$. So that,
if $r\ge2,s>1$ then H\"older inequality, the boundedness of $H_k$ and the above
inequality yield \beqno{Hkey}
\left(\itQ{\qqq}{|H_k|^r}\right)^{\frac1s} \le
C(M)\left(\itQ{\qqq}{|H_k|^2}\right)^{\frac{1}{s}} \le
C(M,r,s,\eg_1)\itQ{\qqq}{|H_k|^2}.\eeq

Since a.4) was not used until we obtain \mref{bbA}, we now need only look at the integral of $\myprod{\bbA,DH_k\fg^2}$
in \mref{bbA}, which reads \beqno{bbAz}\sup_{s\in[-(\frac23)^p,0]}
\iidx{B_\frac23}{H_k^2\fg^2} +\itQ{\qqq}{\myprod{\bbA,DH_k\fg^2}}
\le F_k,\eeq where $\bbA=A_k(u_k,Du_k)-A_k(U_k,DU_k)$, and
$F_k=F_k(h,r,\rg)$ which is defined by \mref{key2} and can be arbitrarily small if $k$ is large
and $h,r,\rg$ small (uniformly in $k$) thanks to the argument in the proof of \refprop{heat-aa} without using a.4).

We now consider the following two cases.

{\bf The case $p>2$:} We write $Q_\frac23 = E_u\cup E_v$ where
$$E_u = \{z\,:\, |Du_k(z)|\le |DU_k(z)|\}, \, E_v = \{z\,:\,
|DU_k(z)|< |Du_k(z)|\}.$$  We also write
$$\itQ{Q_\frac23}{\myprod{\bbA,DH_k\fg^2}} =
\itQ{E_u}{\myprod{\bbA,DH_k\fg^2}}+\itQ{E_v}{\myprod{\bbA,DH_k\fg^2}}.$$
On $E_u$, we have $$\myprod{\bbA,DH_k} =\myprod{
A_k(U_k,DU_k)-A_k(U_k,Du_k)
+A_k(U_k,Du_k)-A_k(u_k,Du_k),D(U_k-u_k)}.$$ By \mref{AAmono},
 as $|Du_k(z)|\le |DU_k(z)|$ and $p>2$, it follows that
\beqno{AAmono1}\myprod{A_k(U_k,DU_k)-A_k(U_k,Du_k), DU_k-Du_k} \ge
\llg_0|Du_k|^{p-2}|DU_k-Du_k|^2,\eeq and this term will be kept on
the left of \mref{bbAz}. On the other hand, as $A_k(u,\zeta)$ is Lipschitz in $u$, 
$$\barr{l}|\myprod{A_k(U_k,Du_k)-A_k(u_k,Du_k),DH_k}| \le
C|H_k||Du_k|^{p-1}|DH_k| \\=
C|H_k||Du_k|^{\frac{p}{2}}|Du_k|^{\frac{p-2}{2}}|DH_k|\le
\frac{\llg_0}{4}|Du_k|^{p-2}|DH_k|^2 + C|H_k|^2||Du_k|^p.\earr$$ The
first term on the right can be absorbed into \mref{AAmono1}.

Interchange the roles of $u_k,U_k$  and apply a similar treatment
for the integral of $\myprod{\bbA,DH_k\fg^2}$ over $E_v$ in the
above argument to see that \mref{bbAz} gives
$$\barr{lll}\lefteqn{\sup_{s\in[-(\frac23)^p,0]}
\iidx{B_\frac23}{H_k^2\fg^2} + \itQ{E_u}{|Du_k|^{p-2}|DH_k|^2\fg^2}
+\itQ{E_v}{|DU_k|^{p-2}|DH_k|^2\fg^2}
}\hspace{3cm}&&\\&\le&C\itQ{E_u}{|Du_k|^p|H_k|^2\fg^2} +
C\itQ{E_u}{|DU_k|^p|H_k|^2\fg^2}+F_k.\earr$$

Since $|Du_k|\le|DU_k|$ on $E_u$, the above yields
$$\sup_{s\in[-(\frac23)^p,0]} \iidx{B_\frac23}{H_k^2\fg^2} \le
C\itQ{\qqq}{|DU_k|^p|H_k|^2\fg^2} +F_k$$

Now, using the higher integrability of $Du_k$ of \reflemm{duLq} we
also have the $L^q$ estimate, with $q>p$, for $DU_k$ on $\qqq$ by
extending $U_k=u_k$ beyond the boundary of $\qqq$. In fact, this comes from
\mref{dukest} (for cylinders centered on the boundary of $\qqq$) to obtain a
reverse H\"older inequality for $DU_k$ and then use the fact that
$Du_k$ is in $L^q$. Thus, by H\"older's inequality and \mref{Hkey}
$$\itQ{\qqq}{|DU_k|^p|H_k|^2\fg^2} \le
C(M)\left(\itQ{\qqq}{|H_k|^r\fg^r}\right)^\frac 1s \le
C(\eg_1,M)\itQ{\qqq^\tau}{|H_k|^2}.$$

Again, it is easy to see that the above argument still holds if
$\qqq$ is replace by $\qqq^\tau=\qqq\cap\{(x,t):t\le\tau\}$ for any
$\tau\le0$. Sending $\rg,r$ to zero to obtain

\beqno{bbAnew}\sup_{s\in[-(\frac23)^p,\tau]} \iidx{B_\frac23}{H_k^2}
\le C(\eg_1,M)\itQ{\qqq^\tau}{|H_k|^2} +F_k.\eeq

Setting $$y_k(\tau)=\sup_{s\in[-(\frac23)^p,\tau]}
\iidx{B_\frac23}{H_k^2}$$ the above \mref{bbAnew} becomes $y_k'\le Cy_k + F_k$
with $y(-(\frac23)^p), F_k$ can be arbitrarily small. By the
Gronwall inequality, we see that $H_k\to 0$ in $L^2(\qqq)$ as well
as in $L^p(\qqq)$ because $p>2$ and $H_k$ is bounded. Our desired contradiction is obtained and the proof is complete for the case $p>2$.

{\bf The case $2n/(n+2)<p<2$:} For each $z\in Q$, let us denote by
$u(z),v(z)$ which are either $u_k$ or $U_k$ such that, with a slight
abuse of notation here, $Du(z)=\max\{Du_k(z),DU_k(z)\}$ and $Dv(z)=\min\{Du_k(z),DU_k(z)\}$. We write
$Q_\frac23 = E_u\cup E_v$ where
$$E_u = \{z\,:\, \frac12|Du(z)|\le |Du(z)-Dv(z)|\}, \, E_v = \{z\,:\,
|Du(z)-Dv(z)|< \frac12|Du(z)|\}.$$ Again, we consider the integral
of $\myprod{\bbA,DH_k\fg^2}$ in \mref{bbAz}. On $E_u$, we write $H=u-v$, $DH=Du-Dv$ and
$$\myprod{\bbA,DH_k} = \myprod{A_k(u,Du)-A_k(u,Dv)
+A_k(u,Dv)-A_k(v,Dv),DH}.$$ Because $p<2$,
$\min\{|Du_k|^{p-2},|Du_k|^{p-2}\}=|Du|^{p-2}$. By \mref{AAmono} and
the fact that  $|Du|^{p-2}\ge 2^{p-2}|Dv-Du|^{p-2}$, we have
\beqno{AAmono1z}\myprod{A_k(u,Du)-A_k(u,Dv), Du-Dv} \ge
C\llg_0|Dv-Du|^p,\eeq and this term will be kept on the left of
\mref{bbAz}. On the other hand, as $A_k(v,\zeta)$ is Lipschitz in $v$ 
$$\barr{lll}|\myprod{A_k(u,Dv)-A_k(v,Dv),DH}| &\le&
C|H||Dv|^{p-1}|DH| \le \frac{C\llg_0}{4}|DH|^p +
C|H|^\frac{p}{p-1}||Dv|^p \\&\le& \frac{C\llg_0}{4}|DH|^p +
C(M)|H|^2||Dv|^p.\earr$$ Here, we have used the fact that $H$ is bounded in the last inequality. The first term on the right can be absorbed
into \mref{AAmono1z}.

On $E_v$ we note that $|Dv-Du|< \frac12|Du|$ implies
$|Dv|\ge\frac12|Du|$ (otherwise, $|Dv-Du|\ge |Du|-|Dv|>\frac12|Du|$
contradicting the definition of $E_v$). We now write $H=v-u$,
$DH=Dv-Du$ and
$$\myprod{\bbA,DH_k} = \myprod{A_k(v,Dv)-A_k(v,Du)
+A_k(v,Du)-A_k(u,Du),DH}.$$ We have by \mref{AAmono}
$$\myprod{A_k(v,Du)-A_k(v,Dv), Dv-Du} \ge
C\llg_0|Du|^{p-2}|Dv-Du|^2.$$ Again, this term will stay on the left
of \mref{bbAz}. Meanwhile, as  $|Dv|\sim|Du|$ on $E_v$, we have
$$\barr{ll}\lefteqn{|\myprod{A_k(v,Du)-A_k(u,Du),DH}|\le
C|H||Du|^{p-1}|DH| \le}\hspace{2cm}&\\
&C|H||Dv|^{\frac{p}{2}}|Du|^{\frac{p-2}{2}}|DH|\le
\frac{C\llg_0}{4}|Du|^{p-2}|DH|^2 + C|H|^2|Dv|^p.\earr$$

Combining the above estimates and noting that $|Dv|\le|DU_k|$, we
derive from \mref{bbAz}
\beqno{bbAnew1z}\barr{lll}\sup_{s\in[-(\frac23)^p,0]}
\iidx{B_\frac23}{H_k^2\fg^2} &+& \itQ{E_u}{|DH_k|^p\fg^2}
+\itQ{E_v}{|Du|^{p-2}|DH_k|^2\fg^2}\\
&\le&C(M)\itQ{\qqq}{|DU_k|^p|H_k|^2\fg^2}+F_k.\earr\eeq

As in the case $p>2$, because $2n/(n+2)<p<2$ the higher
integrability of $DU_k$ is available, we can derive a Gronwall
inequality and see that $H_k\to 0$ in $L^2(\qqq)$. Since $H_k$ is bounded and $DU_k$ is $L^q(\qqq)$, an application of H\"older's inequality shows that the right hand side of the above inequality tends to zero as $k\to\infty$. Hence,

\beqno{bbAnew1z5}\itQ{E_u}{|DH_k|^p\fg^2}
+\itQ{E_v}{|Du|^{p-2}|DH_k|^2\fg^2}\to 0 \quad \mbox{as } k\to\infty.\eeq

Concerning the
$L^p$ norm of $DH_k$, we observe that
$$\barr{lll}\itQ{E_v}{|DH_k|^p\fg^2} &=&
\itQ{E_v}{|Du|^{\frac{(p-2)p}{2}}|DH_k|^p|Du|^{\frac{(2-p)p}{2}}\fg^2}\\
&\le& \left(\itQ{E_v}{|Du|^{p-2}|DH_k|^2\fg^2}\right)^\frac p2
\left(\itQ{E_v}{|Du|^p\fg^2}\right)^\frac{2-p}{2}.\earr$$ Because
the integral of $|Du|^p=\max\{|Du_k|^p,|DU_k|^p\}$ over $\qqq$ is
bounded, the above and \mref{bbAnew1z5} show that the integral of
$|DH_k|^p$ over $E_u$ and $E_v$, and therefore $\qqq$, tends to
zero. By Sobolev's imbedding inequality, we see that $H_k\to 0$ in
$L^p(\qqq)$. Our proof is then complete. \eproof

The first alternative \mref{hb3a} in the above propositions is not as
useful as \mref{ha3a} for  our later  proof of the H\"older continuity for weak solutions. To this end, we will show that \mref{hb3a} allows us to approximate the considered $u$ by solutions to those of
nice systems whose coefficients $A(v,Dv)$ do not involve with the
solutions $v$. The proof of this fact for the singular case ($p<2$) is much more involved and will be reported in a forthcoming work. Here, we present only the result when $p>2$.

\bprop{heat-aab}Assume that $p>2$ and a.1)-a.3) and a.4) or a.4')
hold. For any given $M,\eg>0$ and $\bg>1$ there exists $\dg\in(0,1]$
that depends only on $\llg,\LLg,\LLg_0,\eg,\bg$ such that if
$A\in\ccA, B\in\ccB$ and
 $u\in V(Q_R)$ satisfying
\beqno{ha1az} \mitQ{Q_R}{|u-u_R|^2} + R^p\mitQ{Q_R}{|Du|^p}\le
M,\eeq and\beqno{ha1aaz} |
\itQ{Q_R}{-u\fg_t+\myprod{A(u,Du),D\fg}-\myprod{B(u,Du),\fg}}|\le\dg
\|Du\|_{L^\frac pq(Q_{R})}\|D\fg\|_{L^p(Q_R)}\eeq for any
$Q_R\subset Q$, $\fg\in V_p^0(Q)$, then there exists $v$ in $V_p(Q_{R/2})$ such that
$\|v\|_{L^\infty(Q_{R/2})}\le C(\llg,\LLg,\llg_0,M)$ and
 \beqno{ha3az}
\left\{\barr{l}\mitQ{Q_{R/2}}{|v-u|^2}+\mitQ{Q_{R/2}}{|v-u|^p} \le
\eg R^{p}\mitQ{Q_R}{|Du|^p},\\
\mitQ{Q_{R/2}}{|Dv|^p}\le \bg\mitQ{Q_{R/2}}{|Du|^p}.\earr\right.\eeq

Moreover, $v$ satisfies
\beqno{ha3aaz}\itQ{Q_{R/2}}{-v\fg_t+\myprod{\widetilde{A}(v,Dv),D\fg}-\myprod{\widetilde{F}(v,Dv),\fg}}=0,\quad
\forall \fg\in C_0^1(Q_{R/2}),\eeq where either
$\widetilde{A}(v,Dv)=A(v,Dv)$ and $\widetilde{B}(v,Dv)=B(v,Dv)$ or
$\widetilde{A}(v,Dv)=A(c,Dv)$ and $\widetilde{B}(v,Dv)=B(c,Dv)$ for
some constant vector $c\in \RR^m$.\eprop

\bproof Again, we will only discuss the case $B\equiv0$ here. As in the proof of \refprop{heat-aac}, by scaling we can assume that $R=1$ and need only consider
the case when
\beqno{duto0}\limsup_{k\to\infty}\mitQ{Q_1}{|Du_k|^p}=0.\eeq

We now look at the solution of
\beqno{vconst} \itQ{Q_{\frac23}}{[-v_k\fg_t +
A_k((u_k)_{Q_1},Dv_k)D\fg]}=0 \quad \forall \fg\in V_p^0(Q_\frac23),\eeq
and $v_k=u_k$ on $S_{\frac23}$. By testing \mref{vconst} with
$v_k-u_k$ we easily see that \beqno{dvk} \mitQ{Q_\frac23}{|Dv_k|^p}\le C \mitQ{Q_1}{|Du_k|^p}.\eeq

As before, by subtracting the equations for $u_k,v_k$ and testing
the result with $(v_k-u_k)$, we have for $H_k=v_k-u_k$ the
following \beqno{uvzzz}\sup_t\iidx{B_\frac23}{H_k^2} +
\itQ{\qqq}{ \myprod{A_k((u_k)_{Q_1},Dv_k)-A_k(u_k,Du_k),DH_k}} \le
0.\eeq We now write
$$\barr{lll}A_k((u_k)_{Q_1},Dv_k)-A_k(u_k,Du_k)
&=&A_k((u_k)_{Q_1},Dv_k)-A_k((u_k)_{Q_1},Du_k)\\&&+A_k((u_k)_{Q_1},Du_k)-A_k(u_k,Du_k)\earr$$
and keep the first difference on the left of \mref{uvzzz}. Using the
ellipticity and Lipschitz property of $A_k$, we obtain
\beqno{uvzzz1}\sup_t\iidx{B_\frac23}{H_k^2} \le C\itQ{\qqq}{
|Du_k|^p|u_k-(u_k)_{Q_1}|^2}.\eeq

We now make use of the $L^q$ estimate for $Du_k$, see \mref{lqdu} of
\reflemm{duLq} and \refrem{duLq-rem}, to find some $q>p$  and
$\sg\ge1$ such that
$$\left(\mitQ{\qqq}{|Du_k|^q}\right)^\frac pq \le
C\left(\mitQ{Q_1}{|Du_k|^p}\right)^\sg.$$ Hence, for $r=(p/q)'$
\beqno{Hkestz}\sup_t\iidx{B_\frac23}{|H_k|^2} \le
C\left(\itQ{Q_1}{|Du_k|^p}\right)^\sg\left(\itQ{\qqq}{|u_k-(u_k)_{Q_1}|^{2r}}\right)^\frac{1}{r}.\eeq

Since $u_k$ is bounded,
$$\left(\itQ{\qqq}{|u_k-(u_k)_{Q_1}|^{2r}}\right)^\frac{1}{r} \le
C(M)\left(\itQ{\qqq}{|u_k-(u_k)_{Q_1}|^{2}}\right)^\frac{1}{r}.$$

By \mref{duto0} and an application of the Poincar\'e inequality, we see
that the above quantities tend to 0 as $k\to\infty$. Sending
$\rg,\tau$, the parameters defining $\fg$, to 0 in \mref{Hkestz} and using the boundedness of the integral of $|Du_k|^p$, we
derive \beqno{Hksmall} \sup_t\iidx{B_\frac23}{|H_k|^2} \le
\eg\itQ{Q_1}{|Du_k|^p}\eeq for any given $\eg>0$ when $k$ is sufficiently large. Integrating the
above in $t$, we obtain $$\itQ{\qqq}{|u_k-v_k|^2}\le
\eg\itQ{Q_1}{|Du_k|^p}.$$

Finally, by applying the interpolation inequality in $x$ to
$(u_k-v_k)$, one has the following
$$\iidx{B_{\frac23}}{|u_k-v_k|^p} \le \eg \iidx{B_{\frac23}}{|Du_k-Dv_k|^p}
+ C(\eg)\mypar{\iidx{B_{\frac23}}{|u_k-v_k|^2}}^{\frac{p}{2}}.$$

Integrating with respect to $t$ and using \mref{dvk} to get
$$\itQ{Q_{\frac23}}{|u_k-v_k|^p} \le
C\eg\itQ{Q_1}{|Du_k|^p} +
C(\eg)\mypar{\sup_t\iidx{B_{\frac23}}{|u_k-v_k|^2}}^{\frac{p}{2}}.$$

By \mref{Hksmall}, we see that $$\itQ{Q_{\frac23}}{|u_k-v_k|^p} \le
C\eg\itQ{Q_1}{|Du_k|^p}.$$

The proof is now complete by rescaling. \eproof

\section{Decay estimates and the proof of \reftheo{regmain}}\label{decaysec}\eqnoset

We will prove in this section that the set $\ccI$ of parameters with which the scaling decay property D) holds for our systems is open. Firstly, we recall the property D): For a given bounded function $v$ in $V_p(Q_{1,1})$  we say that $v$
satisfies a {\em scaling decay property} if the following holds

\bdes \item[D)] Let $M=\sup_{Q_{1,1}} |v|$. For any $R_0>0$ and
$\eta\in(0,1)$, there are
 positive numbers $A,K,L,\ag_0,\og_0$ depending on
$M,\eta$ (with $K,A$ sufficiently large) such that we can define the
following sequences\beqno{scale} R_k=\frac{R_0}{K^k}, \quad
\og_{k+1}=\max\{\eta\og_k, LR_{n}^{\ag_0}\},\quad
S_k=\frac{\og_k}{A}, \quad Q_k=B_{R_k}\times
[-S_k^{2-p}R_k^p,0]\eeq such that $Q_k\subset Q_{1,1}$ for any integer $k\ge0$  and 
\beqno{va} \og_k^p \ge \mitQ{Q_k}{|v-(v)_k|^p},\quad
(v)_k=\mitQ{Q_k}{v}.\eeq \edes

We then consider a family of systems of the form ($\tau\in[0,1]$)
\beqno{sysfam}\left\{\barr{ll}v_t = \Div(A(\tau, v,Dv)) & \mbox{ in $Q_{1,1}$},\\
v=g & \mbox{ on $S_{1,1}$.}\earr\right.\eeq

We defined $\ccI$ to be the collection of $\tau\in[0,1]$ such
that every bounded weak solutions of the above system verifies D). \reftheo{regmain} asserted that
$\ccI$ is open and bounded weak solutions to
\mref{sysfam} with $\tau\in \ccI$ are H\"older continuous. Its proof goes as follows.

{\bf Proof of \reftheo{regmain}:} Fix a $\mu\in \ccI$. We will show that if $|\nu-\mu|$ is
sufficiently small then $\nu\in \ccI$. That is, every bounded weak
solution $u$ of \mref{sysfam} with $\tau=\nu$ will satisfy D). Now, let $u$ be such a solution
and $M=\sup_Q|u|$.

The new set of parameters $A,K,L,\ag_0,R_0,\{\og_k\}$ in D) for $u$  will be
determined in the course of our calculation and depend from that of the reference system \mref{sysfam} when
$\tau=\mu$. 

By induction, let us start with a positive $\og_k$, says
$k=0$, such that

\beqno{ua} \og_k^p \ge \mitQ{Q_k}{|u-(u)_k|^p}.\eeq

In the sequel, for any $t>0$ we will denote by $tQ_k$ the cylinder
with radius $t R_k$ and concentric with $Q_k$. From the Caccioppoli
inequality, see \reflemm{CACCineq} with $\rg=S_k^{2-p}R_k^p$,
$$R_k^p\mitQ{\frac12 Q_{k}}{|Du|^p} \le C
\mitQ{Q_{k}}{|u-(u)_k|^p}+C\mitQ{Q_{k}}{S_k^{p-2}|u-(u)_k|^2}.$$  An application of Young's inequality to the integrand in the second
term on the right and \mref{ua} we can find a constant $C_1$
such that \beqno{dua}R_k^p\mitQ{\frac12Q_k}{|Du|^p} \le
C\mitQ{Q_k}{|u-(u)_k|^p} + CS_k^p \le
C_1(1+\frac{1}{A^p})\og_k^p.\eeq

For any given $\eg>0$, if $|\mu-\nu|$ is sufficiently small
(depending only on $M,\eg$) we apply the approximation result,
\refprop{heat-aab}, in $\frac12Q_k$ to obtain a "nice" solution $v$
satisfying \mref{sysfam} with $\tau=\mu$ or a similar system with $v$ being replaced by a constant vector such that $\sup_{\frac14
Q_k}|v|\le C(M)$ and (in combination with \mref{ua} and \mref{dua})
$$
\barr{lll}\mitQ{\frac14Q_k}{|v-(v)_k|^p} &\le&\mitQ{\frac14Q_k}{|v-(u)_k|^p} \\&\le&
\mitQ{\frac14Q_k}{|v-u|^p} + \mitQ{\frac14Q_k}{|u-(u)_k|^p} \\&\le&
\eg R_k^p\mitQ{\frac12Q_k}{|Du|^p} + C_2\og_k^p \le
C_3(1+\frac{1}{A^p})\og_k^p\le C_4\og_k^p.\earr$$ Here, $C_4=2C_3$ if
$A\ge1$. We then take $\hat{\og}_k =C_4\og_k$ ($k=0$) and apply the
assumption D) on any solution $v$ of \mref{sysfam} with $\hat{\eta}$
sufficiently small to find $\hat{A}, \hat{K},\hat{L}$ depending on
$\sup_Q|v|$, and therefore $M$, and construct the sequence
$\{\hat{\og}_k\}$ such that the relations in \mref{scale} hold and

$$\mitQ{\hat{Q}_{k}}{|v-(v)_{k}|^p} \le \hat{\og}^p_{k} \quad \forall k.$$

The new constants for $u$ will be chosen such that $K=\hat{K},
C_4A=\hat{A}$, with $\hat{A}$ being large and $A\ge1$. The constants
$L,\ag_0$ will be determined later (using the constant $C_4,A$) so
that $Q_{k+1}\subset\hat{Q}_{k+1}$ and $\og_{k+1}=\max\{\eta\og_k,
LR_{k}^{\ag_0}\}$.

Choosing $K$ large (or equivalently $\hat{K}$) depending on $\eta$ such that,
as $\og_{k+1}\ge \eta\og_k$, we have
$$\frac{4^p}{K^p} \og_k^{p-2} \le \og_{k+1}^{p-2} \Leftrightarrow S_{k+1}^{2-p}\left(\frac{R}{K^{n+1}}\right)^p \le
S_k^{2-p}\left(\frac{R}{4K^{n}}\right)^p\Rightarrow
Q_{k+1}\subset \frac14Q_k.$$

Noting that $\frac14Q_k$ is scaled by $R_k^p$ in the $t$ direction,
we then deduce
$$\barr{lll}\mitQ{Q_{k+1}}{|u-v|^p} &\le&
\frac{\og_{k+1}^{p-2}}{\og_k^{p-2}}(2K)^{n+p}\mitQ{\frac14Q_{k}}{|u-v|^p}
\le \frac{\og_{k+1}^{p-2}}{\og_k^{p-2}}(2K)^{n+p}\eg
R_k^p\mitQ{\frac12Q_k}{|Du|^p}
\\&\le& C_5\eg (2K)^{n+p}\og_{k+1}^{p-2}\og_k^2.\earr$$
Here, we have used \mref{dua}. Since $Q_{k+1}\subset\hat{Q}_{k+1}$, it follows that

$$\barr{lll}\mitQ{Q_{k+1}}{|u-(u)_{k+1}|^p} &\le& C_6\mitQ{Q_{k+1}}{|u-v|^p} +
C_6\mitQ{Q_{k+1}}{|v-(v)_{k+1}|^p}\\& \le&
C_6\mitQ{Q_{k+1}}{|u-v|^p} +
C_6\frac{|\hat{Q}_{k+1}|}{|Q_{k+1}|}\mitQ{\hat{Q}_{k+1}}{|v-(v)_{k+1}|^p}.\earr$$

The first term on the right is estimated as follows. We choose $\eg$
small, depending on $K$ and thus $M$, such that for any given
$\eg'>0$ we have via Young's inequality
$$C_5\eg (2K)^{n+p}\og_{k+1}^{p-2}\og_k^2 \le \frac12\og_{k+1}^p +
C_7(\eg K^{n+p})^{p/2}\og_k^p \le \frac12\og_{k+1}^p +
\eg'\og_k^p.$$

Meanwhile, since $v$ verifies \mref{va} in D)
$$C_6\frac{|\hat{Q}_{k+1}|}{|Q_{k+1}|}\mitQ{\hat{Q}_{k+1}}{|v-(v)_{k+1}|^p}
=
C_6\frac{S^{p-2}_{k+1}}{\hat{S}^{p-2}_{k+1}}\mitQ{\hat{Q}_{k+1}}{|v-(v)_{k+1}|^p}
\le C_8\og_{k+1}^{p-2}\hat{\og}_{k+1}^2 \le \frac14\og_{k+1}^p$$ if
 $\sqrt{C_8}\hat{\og}_{k+1} \le \og_{k+1}$. Here, from the definition of $C_4A=\hat{A}$, $C_8=C_6C_4^{2-p}$ which depends only on $n$.

Recall that we also require $Q_{k+1}\subset \hat{Q}_{k+1}$. To this
end, as $K=\hat{K}$, we need $S_{k+1}=\frac{\og_{k+1}}{A}\ge
\frac{\hat{\og}_{k+1}}{\hat{A}}=\hat{S}_{k+1}$ or $\og_{k+1}\ge
C_4\hat{\og}_{k+1}$. This and the requirement $ \le \og_{k+1} \ge \sqrt{C_8}\hat{\og}_{k+1}$ are possible by choosing $L$ sufficiently
large or $R_0$ is small with
 $\ag_0<\hat{\ag}_0$ and $\hat{\eta}$ small
so that $$\og_{k+1} = \max\{\eta\og_k,LR^{\ag_0}\} \ge
\max\{\sqrt{C_8},C_4\}\max\{\hat{\eta}\hat{\og}_k,\hat{L}R^{\hat{\ag}_0}\}=\max\{\sqrt{C_8},C_4\}\hat{\og}_{k+1}.$$
We should note that once this requirement is fulfilled for $k=0$
then $\hat{\eta}$ is fixed and the above relation holds for all
$k\ge1$. By induction, we then define the sequence $\{\og_k\}$ for
$u$.

Hence,
$$\mitQ{Q_{k+1}}{|u-(u)_{k+1}|^p} \le \frac12\og_{k+1}^p +
\eg'\og_k^p+ \frac{1}{4}\og^p_{k+1} \le \og_{k+1}^p \quad \mbox{ for
all $k\ge0$}.$$

This shows that $u$ satisfies the same properties D) and $\ccI$ is
open. In addition, an algebraic argument similar to \cite[Proposition 3.1,
p.44]{Dib} applies to the sequence $\{\og_k\}$ to get
$$\mitQ{Q_{k}}{|u-(u)_{k}|^p} \le \og_k^p\le CR_k^{p\ag_0}.$$
Moreover, since $\og_k\le\og_0$ and $p>2$, we have
$Q_{R_k,S_0^{2-p}R_k^p}\subset Q_k$ and the above gives
\beqno{Campkey}\mitQ{Q_{R_k,S_0^{2-p}R_k^p}}{|u-(u)_{k}|^p} \le
CR_k^{\ag_0}.\eeq By Campanato imbedding, the above implies that $u$ is H\"older
continuous.  Our proof is then complete. \eproof

\section{$\ccI$ is closed}\eqnoset\label{iclosed}

In this section, we will show that the set $\ccI$ is closed in
$[0,1]$.  To proceed, we take a sequence $\{\nu_k\}$ in $\ccI$ such
that $\nu_k\to\mu$ and show that $\mu\in \ccI$. Thus, let us
consider a bounded weak solution to \mref{fame1a} with $\nu=\mu$. By
II), there is a sequence of H\"older continuous weak  solutions $v_k$ to \beqno{s1dv}v_t = \Div(A(\nu_k,v,Dv))\mbox{ in $Q_1$}, \eeq such that $Dv_k$ converges weakly in
$L^1(Q_1)$ to $Du$. Moreover, the
 $L^\infty$ norms of $v_k$'s are bounded
uniformly in terms of that of $u$, and  We will derive uniform estimates for various
integral norms of $Dv_k$ in terms of the $L^\infty$ norm of $u$. Once
this is established, we obtain estimates for the derivatives of the
limiting $u$ and its H\"older continuity to conclude that $\mu$
is in $\ccI$.

Let $v$ be any bounded weak solution to \mref{s1dv}. We recall our assumptions here.

Let
$\llg_{\nu,v},\LLg_{\nu,v}$ be the ellipticity constants  for the
matrix $(A^{ij}_{kl})=\frac{\partial A}{\partial\xi}(\nu,v,\xi)$,
that is \beqno{A-ell}\sum_{i,j=1}^m\sum_{k,l=1}^n
A^{ij}_{kl}\eta^i_k\eta^j_l \ge \llg_{\nu,v}|\eta|^2, \quad
\sum_{i,k}(\sum_{j,l}A^{ij}_{kl}\eta^j_l)^2\le
\LLg_{\nu,v}^2|\eta|^2\eeq for any $\eta\in\RR^{mn}$. Moreover, for some positive constants
$\llg_\nu,\LLg_\nu$ we have
$$\llg_{\nu,v}\ge \llg_\nu|Dv|^{p-2}, \quad \LLg_{\nu,v}\le
\LLg_\nu|Dv|^{p-2}.$$If $n>2$, we also assume that
\beqno{lambda-cond}
\frac{\LLg_{\nu,v}}{\llg_{\nu,v}}<\frac{n}{n-2}.\eeq

We also assume that
there exists a positive constant $a_{\nu,v}$ such that \beqno{MM}
|\frac{\partial A} {\partial v} (\nu,v,\xi)| \le
a_{\nu,v}|\xi|^{p-1} \mbox{ with } 2a_{\nu,v}
M_{\nu,v}(p+n-1)<\sg_0\widehat{\llg}_{\nu}, \eeq where $M_{\nu,v}=\sup_{\qqq}|v|$ and $\sg_0$ is a
fixed number in $(0,1)$ and
$$\widehat{\llg}_{\nu} = (1-\dg^2)\llg_{\nu} \mbox{ and }
\dg=\frac{n-2}{n}\sup\{\frac{\LLg_{\nu,v}}{\llg_{\nu,v}}\,:\,\mbox{
$v$ is a bounded solution}\}.$$
Note that $\widehat{\llg}>0$ thanks to \mref{lambda-cond}. 

Fixing $\nu$ in $I$
and a solution $v$ to \mref{s1dv}, we will denote
$a(v,\zeta)=A(\nu,v,\zeta)$ and also omit the parameter $\nu$ in the
subscripts for $\llg_{\nu,v},\LLg_{\nu,v},a_{\nu,v}$ in the sequel.

The proof of \reftheo{dvp1} relies mainly on the following two
lemmas which establish uniform bounds for the $L^q$ norms of $Dv$.
First of all, we need the following simple consequence of Sobolev's
inequality. For any $q,r>0$, assuming $n>2$ as the case $n=2$ is easy, we have $q +r\frac 2n =
\frac12q\frac{2n}{n-2}\frac{n-2}{n} + r\frac2n$ and by H\"older and
Sobolev's inequalities the following
$$\barr{lll}\io{|V|^{q+r\frac2n}\fg^{2+\frac4n}} &\le&
\left(\io{[|V|^{\frac12q}\fg]^{\frac{2n}{n-2}}}\right)^\frac{n-2}{n}
\left(\io{|V|^{r}\fg^2}\right)^\frac{2}{n}\\&\le&
C\io{|D(|V|^{\frac12q}\fg)|^2}\left(\io{|V|^{r}\fg^2}\right)^\frac{2}{n}.\earr
$$ Therefore, by integrating in $t$\beqno{qsobolev} \itQ{Q_1}{|V|^{q+r\frac2n}\fg^{2+\frac4n}}
\le C\sup_{t\in(-1,0)}
\left(\io{|V|^{r}\fg^2}\right)^\frac{2}{n}\itQ{Q_1}{|D(|V|^{\frac12q}\fg)|^2}.\eeq

In the sequel, we will make use of difference quotients. For any vector valued function $f$, $i=1,\ldots,n$ and real number $h\ne0$, we denote $$ \dg_h^{(i)}f(x,t)=\frac{1}{h}(f(x+he_i,t)-f(x,t)), \quad \mbox{$e_i$ is the unit vector in the $i^{th}$ direction of $\RR^n$}.$$ If an argument holds for any $i$, we will simply omit the superscript $(i)$ in the above notation.

For any $v$ being a weak solution to a nice system, $v$ is H\"older
continuous and the difference $\dg_hv$ weakly solves \beqno{s1ah} (\dg_hv)_t =
\Div(\dg_ha(v,Dv)).\eeq 

We first have the following estimate for "nice" solutions.
\blemm{sl1} Let $v$ be a H\"older continuous weak solution to \mref{s1dv}. For any $\fg\in C^1_0(Q_\frac34)$ there exists a constant $C$ depending on $M=\sup_{\qqq}|v|$ such that \beqno{s2c}
\sup_{t\in(-1,0)}\io{|Dv|^2\fg^2} + \llg \itQ{Q_1}{|Dv|^{p-2}|D^2v|^2\fg^2}
\mbox{ and } \itQ{Q_1}{|Dv|^{p+\frac4n}\fg^{2+\frac4n}}\le C.\eeq
\elemm

\bproof Let $\fg$ be in $C^1_0(\qqq)$. For any function $f$ in $(x,t)$,  $h\ne0$ and $e=e_i$ ($i=1,\ldots,n$), we will write $\dg_h^+f(x,t)=(f(x+he,t)-f(x,t))/h$ and $\dg_h^-f(x,t)=(f(x,t)-f(x+he,t))/h$.
Testing \mref{s1ah} with $\dg_h^+v\fg^2$ and integrating by
part in $x$, we get 
\beqno{diffH}\sup_{t\in(-1,0)}\iidx{\Og_1}{|\dg_h^+v|^2\fg^2}  + \itQ{Q_1}{\myprod{\dg^+_ha(v,Dv),D(\dg_h^+v\fg^2)}} \le \itQ{Q_1}{|\dg_h^+v|^2\fg_t}.\eeq

We then set $$E_0=\{(x,t)\in Q_1\,:\,|Dv(x,t)|\le|Dv(x+h,t)|\}, \quad E_1=Q_1\setminus E_0.$$ 

We now split the integral of $\myprod{\dg_h^+a(v,Dv),D(\dg_h^+v\fg^2)}$ on $Q_1$ into those on $E_0,E_1$. On $E_1$, we have $$\myprod{\dg_h^+a(v,Dv),D(\dg_h^+v\fg^2)}=\myprod{\dg_h^-a(v,Dv),D(\dg_h^-v\fg^2)}=\myprod{\dg_h^-a(v,Dv),D(\dg_h^-v)\fg^2+\dg_h^-vD(\fg^2)}.$$

Concerning the term $\dg_h^-a(v,Dv)$, we write $$\barr{lll}\dg_h^-a(v,Dv)&=&
\frac1h[a(v(x,t),Dv(x,t))-a(v(x,t),Dv(x+h,t))]+\\&&\frac1h[a(v(x,t),Dv(x+h,t))-a(v(x+h,t),Dv(x+h,t))]\\
&=&
\dspl{\int_0^1 \frac{\partial a}{\partial \xi}(v,sDv(x,t)+(1-s)Dv(x+h,t))D\dg_h^-vds}+\\&& 
\dspl{\int_0^1 \frac{\partial a}{\partial v}(sv(x,t)+(1-s)v(x+h,t),Dv(x+h,t))\dg_h^-vds}.\earr$$

Using the fact that $|sDv(x,t)+(1-s)Dv(x+h,t)|\ge |Dv(x+h,t)|$ on $E_1$ and the ellipticity condition of $\partial a/\partial\xi$ we get
$$\myprod{ \frac{\partial a}{\partial \xi}(v,sDv(x,t)+(1-s)Dv(x+h,t))D\dg_h^-vds,D\dg_h^-v}\ge \llg|Dv(x+h,t)|^{p-2}|D\dg_h^-v|^2.$$ This term will stay on the left of \mref{diffH}. On the other hand, bty \mref{MM}, we have
$$\left|\frac{\partial a}{\partial v}(sv(x,t)+(1-s)v(x+h,t),Dv(x+h,t))\dg_h^-v\right|\le |a_v||Dv(x+h,t)|^{p-1}|\dg_h^-v|.$$

Thus, by Young's inequality, we have for any positive $\eg$ the following
$$\barr{lll}\lefteqn{\left|\myprod{\frac{\partial a}{\partial v}(sv(x,t)+(1-s)v(x+h,t),Dv(x+h,t))\dg_h^-v,D(\dg_h^-v)}\right|\le}\hspace{2cm}&&\\&\le&\eg|Dv(x+h,t)|^{p-2}|D(\dg_h^-v)|^2+C(\eg)|a_v|^2|Dv(x+h,t)|^{p}|\dg_h^-v|^2.\earr$$ The above term will be on the right hand side of \mref{diffH}.

Similar argument will apply to the set $E_0$. We then choose $\eg$ sufficiently small and derive from \mref{diffH} and the above estimates
\beqno{s2ah}
\barr{ll}\lefteqn{\sup_{t\in(-1,0)}\iidx{\Og_1}{|\dg_hv|^2\fg^2} +\frac\llg2
\itQ{Q_1}{|v_d|^{p-2}(|D(\dg_h^+v)|^2+|D\dg_h^-v|^2)\fg^2}\le}\hspace{2cm}&\\&
\itQ{Q_1}{\frac{|a_v|^2}{\llg}|v_D|^{p}|V_h|^2\fg^2 +
|v_D|^2(|\fg_t|+\LLg|v_D|^{p-2}|D\fg|^2 +1)},\earr\eeq where $v_D=\max\{|Dv(x,t)|,|Dv(x+h,t)|\}$, $v_d=\min\{|Dv(x,t)|,|Dv(x+h,t)|\}$ and $V_h=\max\{|\dg_h^+v|,|\dg_h^-v|\}$.

Sending $h$ to zero, we get\beqno{s2a}
\barr{ll}\lefteqn{\sup_{t\in(-1,0)}\iidx{\Og_1}{|Dv|^2\fg^2} +\frac\llg2
\itQ{Q_1}{|Dv|^{p-2}|D^2v|^2\fg^2}\le}\hspace{2cm}&\\&
\itQ{Q_1}{\frac{|a_v|^2}{\llg}|Dv|^{p+2}\fg^2 +
|Dv|^2(|\fg_t|+\LLg|Dv|^{p-2}|D\fg|^2 +1)}.\earr\eeq Of course the above argument is justified if $Dv\in L^{p+2}_{loc}$. This fact will be proven  in \reflemm{DvLqloc} following this proof.

We now estimate the integral of $|Dv|^{p+2}$ in \mref{s2a}. By
integrating by parts in $x$
$$\itQ{Q_1}{|Dv|^{p+2}\fg^2} = \itQ{Q_1}{vD(|Dv|^{p+1}\fg^2)} \le M\itQ{Q_1}{(|D^2v||Dv|^p\fg^2+
|Dv|^{p+1}\fg|D\fg|)},$$ where $M=M_{\nu,v}=\sup_{Q_\frac34}|v|$.

Young's inequality applying to the right then gives
$$\itQ{Q_1}{|Dv|^{p+2}\fg^2} \le \eg\itQ{Q_1}{|Dv|^{p+2}\fg^2} +
C(\eg)M^2\itQ{Q_1}{(|D^2v|^2|Dv|^{p-2}\fg^2 + |Dv|^{p}|D\fg|^2)}.$$
Thus, for $\eg=\frac12$, we obtain
$$\frac{a_v^2}{\llg}\itQ{Q_1}{|Dv|^{p+2}\fg^2} \le
4\frac{a_v^2M^2}{\llg}\itQ{Q_1}{|Dv|^{p-2}|D^2v|^2\fg^2 +
|Dv|^p|D\fg|^2}.$$

Using this in \mref{s2a} and the assumption on the smallness of
$a_vM$ in \mref{MM}, we obtain \beqno{s2b}
\barr{ll}\lefteqn{\sup_{t\in(-1,0)}\io{|Dv|^2\fg^2} + \llg
\itQ{Q_1}{|Dv|^{p-2}|D^2v|^2\fg^2}\le}\hspace{2cm}&\\&
C\itQ{Q_1}{|Dv|^2(|\fg_t|+(|Dv|^{p-2}\LLg+\llg)|D\fg|^2
+1)}.\earr\eeq

By Caccioppoli's inequality we note that \beqno{s3}
\itQ{Q_1}{|Dv|^2(|\fg_t|+|Dv|^{p-2}|D\fg|^2)} \le
C(|\fg_t|,|D\fg|)\itQ{Q_1}{|v|^2} \le C(M).\eeq

The above and \mref{s2b} then imply
$$\sup_{t\in(-1,0)}\iidx{\Og_1}{|Dv|^2\fg^2} +
\itQ{Q_1}{|Dv|^{p-2}|D^2v|^2\fg^2}\le C(M).$$ We now make
use of \mref{qsobolev}, with $V=Dv$ and $q=p,r=2$, and combine with
the above to get the second estimate in \mref{s2c} and complete the
proof. \eproof

To justify the calculation leading to \mref{s2a}, we now show that

\blemm{DvLqloc} Spatial derivatives of H\"older continuous weak  solutions $v$ to \mref{s1dv} are in $ L^{p+2}_{loc}$. \elemm

\bproof Let $v$ be a H\"older continuous weak  solutions to \mref{s1dv}. We will show that at almost every point $z_0=(x_0,t_0) \in Q_1$ with $Dv(z_0)\ne0$ and $R$ is sufficiently small then there is a constant $C$ such that\beqno{dvz1} \itQ{Q_R}{|\dg_{h_k}v|^{p+2}} \le C R^{-2}\itQ{Q_{4R}}{|Dv|^p} \quad \mbox{for some sequence $h_k\to0$}.\eeq  Here, $Q_R= B_R(x_0)\times(t_0-R^2,t_0)$ and $Q_R\subset Q_{4R}\subset Q_1$.
If this is not true then there will be  a sequence  $R_k\to0$   such that $$ \itQ{Q_{R_k}}{|\dg_{h_k}v|^{p+2}} > k R_k^{-2}\itQ{Q_{4R_k}}{|Dv|^p} \quad \mbox{for any sequence $h_k\to0$}.$$ We then choose $h_k=hR_k$  for some sufficiently small and positive $h$ to be determined later. By scaling, with $v_k(X,T)=v(x_0+R_kX,t_0+R_k^2T)$, we get a sequence of functions $v_k$ on $Q_1$ such that \beqno{dvz2} \itQ{Q_\frac14}{|\dg_{h}v_k|^{p+2}} \ge k\itQ{Q_1}{|Dv_k|^p}.\eeq 

Since $v$ is H\"older continuous and $R_k\to0$ we see that $v_k$ can be arbitrarily close to $(v_k)_{Q_1}$ on $Q_1$. Thus, $v_k$ approximately solves the following system \beqno{dvz3} U_t = \Div(a( (v_k)_{Q_1}, DU)).\eeq For such system, which does not explicitly depend on $U$, we can find (see \refrem{UHest} at the end of this section) a function $C(x)$ which is bounded if $x$ is bounded such that
\beqno{UH}\itQ{Q_\frac14}{|DU|^{p+2}} \le C\left(\itQ{Q_\frac12}{|DU|^p}\right).\eeq 

For sufficiently small $h>0$, the above yields
$$\itQ{Q_\frac14}{|\dg_{h}U|^{p+2}} \le C\left(\itQ{Q_\frac12}{|DU|^p}\right).$$

Our approximation results then give a sequence $\{U_k\}$ of weak solutions to \mref{dvz3} satisfying $$\itQ{Q_\frac12}{|DU_k|^p}\le c\itQ{Q_\frac12}{|Dv_k|^p},$$  and $U_k-v_k\to0$ in $L^{p}(Q_\frac14)$ (as well as in $L^{q}(Q_\frac14)$ for any $q>1$ because $U_k,v_k$ are bounded). Moreover, $DU_k-Dv_k\to0$ weakly in $L^p(Q_\frac14)$.

Together, we have the following estimates $$k\itQ{Q_1}{|Dv_k|^p} \le \itQ{Q_\frac14}{|\dg_{h}v_k|^{p+2}} , \quad \itQ{Q_\frac14}{|\dg_{h}U_k|^{p+2}} \le C\left(\itQ{Q_\frac12}{|Dv_k|^p}\right).$$
With $h$ being fixed and $U_k,v_k$ being bounded, we have $\dg_{h}(U_k)-\dg_{h}(v_k)\to0$ in $L^{q}$ for all $q>1$ as  $k\to\infty$. Moreover, since
$\itQ{Q_\frac12}{|Dv_k|^p}\to Dv(z_0)\ne0$, the above gives a contradiction when $k\to \infty$. Thus, \mref{dvz1} holds almost everywhere on the set where $Dv\ne0$. Finally, by sending $h$ to 0 it is easy to see that \mref{dvz1} implies $Dv\in L^{p+2}_{loc}$.
\eproof

To get estimates for higher powers of $|Dv|$, we need the following
lemma.

\blemm{sl2} Let $v$ be a H\"older continuous weak solution to \mref{s1dv} and $\ag$ be a positive number. Assume that
$\frac{\ag}{2+\ag}=\dg_{\ag,v} \frac{\llg_v}{\LLg_v}$ for some
$\dg_{\ag,v}\in(0,1)$ and \beqno{small} 2a_v
M(p+\ag+1)<\sg_0\widehat{\llg}_v, \quad \mbox{ with }
\widehat{\llg}_v=(1-\dg_{\ag,v}^2)\llg_v,\quad \sg_0\in(0,1). \eeq
If $\fg\in C_0^1(Q_\frac34)$ then
\beqno{main}\itQ{Q_1}{|Dv|^{p+\ag+(2+\ag)\frac2n}\fg^{2+\frac4n}}
\le C\left(\itQ{Q_1}{|Dv|^{2+\ag}(|\fg_t|+|Dv|^{p-2}|D\fg|^2
+1)}\right).\eeq\elemm

\bproof To proceed, we recall the following facts in \cite{ka}. From
the ellipticity condition of $\frac{\partial a}{\partial
\xi}=(A^{ij}_{kl})$ we have for $\kg_v=\llg_v/\LLg_v^2$ and
$\nu_v=\llg_v/\LLg_v$ that
$$\sum_{i,k}(\eta^k_i-\kg A^{ij}_{kl}\eta^j_l)^2 \le
(1-2\kg_v\llg_v+\kg_v^2\LLg_v^2)|\eta|^2=(1-\nu_v^2)|\eta|^2.$$

The lemma \cite[p.677]{ka} then gives $D\zeta
D(\zeta|\zeta|^\ag)\ge\mu^\frac12(\ag)|D\zeta||D(\zeta|\zeta|^\ag)|$
for any $\zeta:\RR^n\to\RR^{nm}$ and $\mu(\ag)=1-(\frac{\ag}{2+\ag})^2$.
Therefore, with $\zeta=Dv$
$$\barr{lll}\sum \frac{\partial a}{\partial
\xi}D^2vD(Dv|Dv|^\ag) &=& \frac{1}{\kg_v}[\kg_v\sum(\frac{\partial a}{\partial
\xi}D^2v-D^2v)D(Dv|Dv|^\ag)] + D(Dv)D(Dv|Dv|^\ag)\\&\ge&
\frac{1}{\kg_v}(\mu^\frac12(\ag)-(1-\nu_v^2)^\frac12)|D(Dv)||D(Dv|Dv|^\ag)|.\earr$$

Thus, if $\frac{\ag}{2+\ag}=\dg_{\ag,v}\frac{\llg_v}{\LLg_v}$ for
some $\dg_{\ag,v}\in(0,1)$ then the constant in the right hand side is
$$\frac{1}{\kg_v}(\mu(\ag)^\frac12-(1-\nu_v^2)^\frac12)=
\frac{1}{\kg_v}\frac{\nu_v^2-\frac{\ag^2}{(2+\ag)^2}}{(\mu(\ag)^\frac12+(1-\nu_v^2)^\frac12)}
\ge
(1-\dg_{\ag,v}^2)\frac{\nu_v^2}{\kg_v}=(1-\dg_{\ag,v}^2)\llg_v.$$
Hence, by the assumption on $\llg_v$, we get
\beqno{D2v}\barr{lll}\sum \frac{\partial a}{\partial
\xi}D^2vD(Dv|Dv|^\ag)&\ge&
\frac1{\kg_v}(\mu^\frac12(\ag)-(1-\nu_v^2)^\frac12)|D(Dv)||D(Dv|Dv|^\ag)\\&\ge&
(1-\dg_{\ag,v}^2)\llg_v|Dv|^{\ag}|D^2v|^2=\widehat{\llg}_v|Dv|^{p-2+\ag}|D^2v|^2.\earr\eeq

The following calculation will be rigorously justified by using difference quotient operator $\dg_h$, as in the previous lemmas, in place of the differentiation $D$ below. However, in order to be more suggestive, we will write \mref{s1ah} formally as 

\beqno{s1a} (Dv)_t =
\Div(\frac{\partial a}{\partial \xi}(v,Dv)D^2v+\frac{\partial
a}{\partial v}(v,Dv)Dv) .\eeq

Testing \mref{s1a} with $Dv|Dv|^\ag\fg^2$ to obtain (compare with \mref{s2ah})
$$\barr{ll}\lefteqn{\sup_t\iidx{\Og_1}{|Dv|^{2+\ag}\fg^2} +
\itQ{Q_1}{\myprod{\frac{\partial a}{\partial
\xi}(v,Dv)D^2v,D(Dv|Dv|^\ag)\fg^2}}\le\itQ{Q_1}{|Dv|^{2+\ag}|\fg_t|}
}\hspace{1cm}&\\& +\itQ{Q_1}{|\myprod{\frac{\partial a}{\partial
\xi}(v,Dv)D^2v,Dv|Dv|^\ag D\fg\fg}|
+\left|\frac{\partial{a}}{\partial
v}\right|(|Dv|^{2+\ag}\fg^2+|Dv|^{1+\ag}|\fg D\fg|)} .\earr$$

Using \mref{D2v} and Young's inequality, we
deduce\beqno{ss2a}
\barr{ll}\lefteqn{\sup_t\iidx{\Og_1}{|Dv|^{2+\ag}\fg^2} +
\widehat{\llg}_v \itQ{Q_1}{|D^2v|^2|Dv|^{p-2+\ag}\fg^2}\le}&\\&
\frac{|a_v|^2}{\widehat{\llg}_v}\itQ{Q_1}{|Dv|^{p+2+\ag}\fg^2} +
C\itQ{Q_1}{|Dv|^{2+\ag}(|\fg_t|+|Dv|^{p-2}|D\fg|^2 +1)}.\earr\eeq

Again, since $v$ is H\"older continous, similar argument as that of \reflemm{DvLqloc} shows that $Dv\in L^{p+2+\ag}_{loc}$ and justifies the above and below calculation. We now estimate the integral of $|Dv|^{p+2+\ag}\fg^2$. By
integrating by parts in $x$, we have
$$\barr{ll}\lefteqn{\itQ{Q_1}{|Dv|^{p+2+\ag}\fg^2} = \itQ{Q_1}{vD(Dv|Dv|^{p+\ag}\fg^2)}}\hspace{2cm}&\\
& \le M\itQ{Q_1}{(p+\ag+1)|D^2v||Dv|^{p+\ag}\fg^2+
|Dv|^{p+1+\ag}\fg|D\fg|},\earr$$
where $M=\sup_{\qqq}|v|$.
Young's inequality applying to the right then gives
$$\barr{ll}\lefteqn{\itQ{Q_1}{|Dv|^{p+2+\ag}\fg^2} \le
(\frac12+\eg)\itQ{Q_1}{|Dv|^{p+2+\ag}\fg^2}+}\hspace{1cm} &\\&
[2M(p+\ag+1)]^2\itQ{Q_1}{|D^2v|^2|Dv|^{p-2+\ag}\fg^2}
+C(\eg)M^2\itQ{Q_1}{ |Dv|^{p+\ag}|D\fg|^2}.\earr$$

We choose $\eg<1/2$ in the above to obtain an estimate for the
integral of $|Dv|^{p+2+\ag}\fg^2$. Using this in \mref{ss2a} and the
assumption \mref{small} on $|a_v M|$ and $\widehat{\llg}_v$, we
obtain \beqno{ss2b}
\barr{ll}\lefteqn{\sup_t\iidx{\Og_1}{|Dv|^{2+\ag}\fg^2} +
(1-\sg_0^2)\widehat{\llg}_v
\itQ{Q_1}{|D^2v|^2|Dv|^{p-2+\ag}\fg^2}\le}\hspace{3cm}&\\&
C\itQ{Q_1}{|Dv|^{2+\ag}(|\fg_t|+\widehat{\llg}|Dv|^{p-2}|D\fg|^2
+1)} .\earr\eeq

The above also gives similar estimate for
$\||Dv|^{(p+\ag)/2}\fg\|_{V(Q_1)}$. Applying \mref{qsobolev}, with
$V=Dv$ and $q=p+\ag,r=2+\ag$, we get the lemma. \eproof

We are now ready to give

{\bf Proof of \reftheo{dvp1}:} By M.2) (see \mref{MM}), we can choose a
number $\bg>n-2$ and some $\dg'\in(0,1)$ such that for $M=\sup_{Q_\frac34}|v|$
$$2a_v M(p+\bg+1)<\sg_0\widehat{\llg}_v, \mbox{ and }
\frac{\bg}{\bg+2}=\dg'\frac{\llg_v}{\LLg_v}>\frac{n-2}{n}.$$
Clearly, starting with $\ag_0\le2$ (thus $p+\ag_0\le p+2$), we can
find finitely many numbers $\ag_0,\ldots,\ag_K$ such that
$$2a_v M(p+\ag_k+1)<\sg_0\widehat{\llg} \quad\mbox{ and } \quad
 \frac{\ag_k}{\ag_k+2}=\dg^{(k)}_{s,v}\frac{\llg_v}{\LLg_v},\quad\dg^{(k)}_{s,v}\le\dg'$$ for any $k\le K$.
Moreover,  it is easy to see that we can also choose
$\dg^{(k)}_{s,v}$ such that $\ag_k\le \ag_{k+1}\le
\ag_k+(2+\ag_k)\frac2n$, $\ag_K=\bg$ and $p+\ag_{k+1}\le
p+\ag_k$.

 Since $Q_\frac57\subset Q_\frac34$, by the Caccioppoli inequality we obtain $\|Du\|_{L^p(Q_\frac57)}\le C(M)$. Using the estimate for $|Dv|^{p+2}$ in
\reflemm{sl1} and a cut-off function $\fg$ for $\qqq,Q_\frac57$ in
\reflemm{sl2}, we see that
$$\itQ{\qqq}{|Dv|^{p+2}} \le C(M).$$

Let the function $\fg$ in \reflemm{sl2} be the cut-off functions for
$Q_{R_k}$ and $Q_{R_{k-1}}$ with $R_k=\frac23-k\frac{1}{6K}$,  $k=1,\ldots,K$. By
induction, the choice of $\ag_k$ allows us to obtain the following estimate from \mref{main} 
\beqno{uni}\itQ{Q_\frac12}{|Dv|^{q_K}} \le C(K,M), \quad
q_K=p+\ag_K+(2+\ag_K)\frac2n.\eeq Note that $q_K>n+p$ because
$\ag_K=\bg>n-2$.

Now, let $u$ be a weak solution which is, by II), approximated by a
sequence $\{v_k\}$ of weak solutions to nice systems and
$Dv_k\rightharpoonup Du$ weakly in $L^1(Q_\frac12)$. By the
semicontinuity of seminorms and \mref{uni}, we have for any
$Q_R\subset Q_\frac12$ that $$ \itQ{Q_R}{|Du|^{q_K}} \le
\liminf_{k\to\infty} \itQ{Q_R}{|Dv_k|^{q_K}} \le C(M).$$ Hence,
with $q=q_K/p$
$$\itQ{Q_R}{|Du|^p} \le
\left(\itQ{Q_R}{|Du|^{q_K}}\right)^\frac1q|Q_R|^{1-\frac1q}\le
C(M)R^{n+p-(n+p)\frac1q}=C(M)R^{n+\ag}.$$ Here,
$\ag=p-(n+p)\frac1q=\frac{p}{q_K}(q_K-n-p)$ is positive. H\"older
continuity for $u$ then follows from the above estimate and the
Poincar\'e inequality in \reflemm{pineq}. Thus, $u$ can be
approximated by solutions to systems that do not explicitly depend
on $u$  and satisfy the property
D). This implies that $u$ satisfies D) and $\ccI$ is closed. \eproof

\brem{UHest} We should note that the estimate \mref{UH} for systems with coefficients independent of their solutions could be derived directly from our proof. Indeed, for such systems $\frac{\partial a}{\partial v}=0$. Therefore, our arguments which lead to \mref{s2a} and \mref{ss2a} would not yield the integral of $|Du|^{p+2+\ag} $ on the right hand sides, and give a similar bound as in \mref{UH} for the integral of $|DU|^{q_K}$ for some $q_K>n+p$. Hence, the calculation in our lemmas is justified.\erem

\bibliographystyle{plain}

\end{document}